\documentclass[11pt]{article}
\usepackage{epic,latexsym,amssymb}
\usepackage{color}
\usepackage{tikz}
\usepackage{amsfonts,epsf,amsmath,leftidx}

\newtheorem{thm}{Theorem}
\newtheorem{lem}{Lemma}

\newtheorem{claim}{Claim}
\newtheorem{subclaim}{Claim}[claim]

\newtheorem{ob}{Observation}

\newcommand{\qed}{$\Box$}

\newcommand{\smallqed}{{\tiny ($\Box$)}}

\newcommand{\cB}{\mathcal{B}}
\newcommand{\cN}{\mathcal{N}}
\newcommand{\cP}{\mathcal{P}}
\newcommand{\cU}{\mathcal{U}}

\newcommand{\cF}{\mathcal{F}}

\newcommand{\cS}{\mathcal{S}}
\newcommand{\cG}{\mathcal{G}}

\newcommand{\barA}{\overline{A}}
\newcommand{\barG}{\overline{G}}
\newcommand{\barT}{\overline{T}}

\newenvironment{unnumbered}[1]{\trivlist \item [\hskip \labelsep {\bf
#1}]\ignorespaces\it}{\endtrivlist}

\newcommand{\proof}{\noindent\textbf{Proof. }}

\let\oldenumerate\enumerate
\renewcommand{\enumerate}{
  \oldenumerate
  \setlength{\itemsep}{0pt}
  \setlength{\parskip}{0pt}
  \setlength{\parsep}{0pt}
}

\begin{document}

\title{Traceability of Connected Domination Critical Graphs}

\author{$^1$Michael A. Henning\thanks{Research supported in part by the University of Johannesburg.}, \, $^2$Nawarat Ananchuen \\ and \\ $^3$Pawaton Kaemawichanurat\thanks{Research
supported by Thailand Research Fund (MRG 6280223)}
\\ \\
$^{1}$Department of Mathematics and Applied Mathematics \\
University of Johannesburg \\
Auckland Park, 2006 South Africa\\
\small \tt Email: mahenning@uj.ac.za \\
\\
$^{2}$Center of Excellence in Mathematics, \\
CHE, Si Ayutthaya Rd., Bangkok 10400, Thailand \\
\small \tt Email: nawarat.ana@mahidol.ac.th \\
\\
$^{3}$Theoretical and Computational Science Center \\
and Department of Mathematics\\
King Mongkut's University of Technology Thonburi \\
Bangkok, Thailand \\
\small \tt Email: pawaton.kae@kmutt.ac.th}

\date{}
\maketitle

\begin{abstract}
A dominating set in a graph $G$ is a set $S$ of vertices of $G$ such that every vertex outside $S$ is adjacent to a vertex in $S$. A connected dominating set in $G$ is a dominating set $S$ such that the subgraph $G[S]$ induced by $S$ is connected. The connected domination number of $G$, $\gamma_c(G)$, is the minimum cardinality of a connected dominating set of $G$. A graph $G$ is said to be $k$-$\gamma_{c}$-critical if the connected domination number $\gamma_{c}(G)$ is equal to $k$ and $\gamma_{c}(G + uv) < k$ for every pair of non-adjacent vertices $u$ and $v$ of $G$. Let $\zeta$ be the number of cut-vertices of $G$. It is known that if $G$ is a $k$-$\gamma_{c}$-critical graph, then $G$ has at most $k - 2$ cut-vertices, that is $\zeta \le k - 2$. In this paper, for $k \ge 4$ and $0 \le \zeta \le k - 2$, we show that every $k$-$\gamma_{c}$-critical graph with $\zeta$ cut-vertices has a hamiltonian path if and only if $k - 3 \le \zeta \le k - 2$.
\end{abstract}


\section{Introduction}

A \emph{dominating set} in a graph $G$ is a set $S$ of vertices of $G$ such that every vertex in $V(G) \setminus S$ is adjacent to at least one vertex in $D$. The \emph{domination number} of $G$, denoted by $\gamma(G)$, is the minimum cardinality of a dominating set of $G$. A graph $G$ is said to be $k$-$\gamma$-critical if $\gamma(G) = k$ and $\gamma(G + uv) < k$ for every pair of non-adjacent vertices $u$ and $v$ of $G$. Such a graph $G$ is called a \emph{domination critical graph}. If $S$ is a dominating set of $G$, we write $S \succ G$, and if $X = \{v\}$, we also write $v \succ G$ rather than $\{v\} \succ G$. The concept of domination and its variations have been widely studied in the literature; a rough estimate says that it occurs in more than 6,000 papers to date. A thorough treatment of the fundamentals of domination theory in graphs can be found in the books~\cite{HHS1,HHS2}.

A \emph{connected dominating set}, abbreviated a CD-set, of a connected graph $G$ is a dominating set $S$ of $G$ such that the subgraph $G[S]$ induced by $S$ is connected. The \emph{connected domination number} of $G$, denoted by $\gamma_c(G)$, is the minimum cardinality of a CD-set of $G$. A CD-set of $G$ of cardinality $\gamma_c(G)$ is called a $\gamma_c$-\emph{set} of $G$. A graph $G$ is said to be $k$-$\gamma_c$-critical if $\gamma_c(G) = k$ and $\gamma_c(G + uv) < k$ for every pair of non-adjacent vertices $u$ and $v$ of $G$. Such a graph $G$ is called a \emph{connected domination critical graph}. If $S$ is a CD-set of $G$, we write $S \succ_c G$, and if $X = \{v\}$, we also write $v \succ_c G$ rather than $\{v\} \succ_c G$. The concept of connected domination was studied at least in the early 1970s, although it was first formally defined by Sampathkumar and Walikar in their 1979 paper~\cite{SaWa79}. Subsequently over the past forty years, the connected domination number has been extensively studied in the literature; a rough estimate says that it occurs in more than 400 papers to date. For a small sample of papers on the connected domination we refer the reader to~\cite{CaWeYu00,DeHaHe13,DuMa09,Ka12,LiZh15,Sa91}.

We remark that the concept of connected domination in graphs is application driven, as evidenced by the earlier papers on the concept. For example, Wu and Li~\cite{WuLi99} show that connected dominating sets are useful in the computation of routing for mobile ad hoc networks. In this application, a minimum connected dominating set is used as a backbone for communications, and vertices that are not in this set communicate by passing messages through neighbors that are in the set.

We also remark that finding connected dominating sets and Steiner trees in a graph are closely related~\cite{CoSt90,DaMo88}.
Moreover, determining the connected domination number of a connected graph $G$ is equivalent to finding the largest possible number of leaves among all spanning trees of $G$.  A \emph{maximum leaf spanning tree} of $G$ is a spanning tree that has the largest possible number of leaves among all spanning trees of $G$, and the \emph{max leaf number}, denoted $\ell_{\max}(G)$, of $G$ is the number of leaves in a maximum leaf spanning tree of $G$. Since $n(G) = \ell_{\max}(G) + \gamma_c(G)$, the problems of a connected dominating set and a maximum leaf spanning tree are closely connected. The maximum leaf spanning tree problem is MAX-SNP hard, implying that no polynomial time approximation scheme is likely~\cite{Ga94}. We remark, however, that both the minimum connected dominating set problem and the maximum leaf spanning tree problem are fixed-parameter tractable~\cite{BiFe14}. The connected dominating set problem is polynomially solvable for distance-hereditary graphs~\cite{DaMo88}.

\subsection{Terminology and Notation}
\label{notation}

For notation and graph theory terminology, we in general follow~\cite{HeYe_book}. Specifically, let $G = (V, E)$ be a graph with vertex set $V=V(G)$ and edge set $E=E(G)$, and let $v$ be a vertex in $V$. A \emph{neighbor} of a vertex is a vertex adjacent to it. The \emph{open neighborhood} of $v$ is the set $N_G(v)$ of all neighbors of $v$, and so $N_G(v)  = \{u \in V \, | \, uv \in E\}$ and the \emph{closed neighborhood of $v$} is $N_G[v] = \{v\} \cup N_G(v)$. A vertex $v$ is said to \emph{dominate} a vertex $u$ in $G$ if $u = v$ or if $u$ is a neighbor of $v$. The \emph{degree} of a vertex $v$ is $|N_G(v)|$ and is denoted by $d_G(v)$. An \emph{end vertex} is a vertex of degree~$1$ and a \emph{support vertex} is a vertex adjacent to an end vertex. For a set $S$ of vertices in $G$, the subgraph induced by $S$ in $G$ is denoted by $G[S]$. If $G$ is a graph, the \emph{complement} of $G$, denoted by $\barG$, is formed by taking the vertex set of $G$ and joining two vertices by an edge whenever they are not joined in $G$.
If the graph $G$ is clear from the context, we omit it in the above expressions. For example, we write $N(v)$ and $N[v]$ rather than $N_G(v)$ and $N_G[v]$, respectively. We use the standard notation $[k] = \{1,\ldots,k\}$.

Two vertices $u$ and $v$ in a graph $G$ are \emph{connected} if there exists a $(u,v)$-path in $G$. A graph $G$ is \emph{connected} if every two vertices in $G$ are connected. We denote the number of components in a graph $G$ by $\omega(G)$. The \emph{distance} $d_G(u,v)$ between two vertices $u$ and $v$ in a connected graph $G$ is the length of a shortest $(u,v)$-path in $G$.
A \emph{hamiltonian cycle} (respectively, \emph{hamiltonian path}) of a graph is a cycle (path) passing through all vertices of the graph. A graph $G$ is \emph{traceable} if it contains a hamiltonian path. Moreover, a graph $G$ is \emph{hamiltonian} if it contains a hamiltonian cycle. For any subgraph $F$ of $G$ and distinct vertices $a$ and $b$ of $G$, $aP_{F}b$ denotes an $(a,b)$-path in $G$ all of whose internal vertices are in $V(F)$. We note that $a$ and $b$ need not be in $V(F)$. If $P$ is an $(a,b)$-path in $G$, we sometimes write the path $P$ by $aPb$ to indicate the start and end vertices of the path $P$.

We denote the \emph{path}, \emph{cycle}, and \emph{complete graph} on $n$ vertices by $P_n$, $C_n$, and $K_n$, respectively, and we denote the \emph{complete bipartite graph} with partite sets of cardinality~$n$ and $m$ by $K_{n,m}$. A \emph{star} is the graph $K_{1,k}$, where $k \ge 1$. The graph $K_{1,3}$ is called a \emph{claw}. A graph $G$ is \emph{claw}-\emph{free} if it does not contain a claw as an induced subgraph. A \emph{tree} is a connected graph with no cycle.

For vertex subsets $X, Y \subseteq V(G)$, we let $N_{Y}(X)$ be the set of all vertices in $Y$ that have a neighbor that belongs to $X$ in $G$, that is, $N_{Y}(X) =  \{y \in Y \mid y \in N_G(x)$  for some $x \in X \}$. For a subgraph $H$ of $G$, we use $N_{Y}(H)$ instead of $N_{Y}(V(H))$ and we use $N_{H}(X)$ instead of $N_{V(H)}(X)$. If $X = \{x\}$, we use $N_{Y}(x)$ instead of $N_{Y}(\{x\})$. The \emph{open neighborhood} of a set $S$ of vertices in $G$ is the set $N_G(S) = \bigcup_{v \in S} N_G(v)$ and its \emph{closed neighborhood} is the set $N_G[S] = N_G(S) \cup S$.

A subset $S \subseteq V(G)$ is a \emph{vertex cut set} of $G$ if the number of components of $G - S$ is more than the number of components of $G$; that is, of $\omega(G - S) > \omega(G)$. In particular, if $S = \{v\}$, then $v$ is called a \emph{cut-vertex} of $G$. We let $\zeta(G)$ be the number of cut-vertices of $G$. When no ambiguity can occur, we write $\zeta$ instead of $\zeta(G)$. A \emph{block} of a graph $G$ is a maximal connected subgraph of $G$ has no cut-vertex of its own. Thus, a block is a maximal $2$-connected subgraph of $G$. Any two blocks of a graph have at most one vertex in common, namely a cut-vertex. A block of $G$ containing exactly one cut-vertex of $G$ is called an \emph{end block}. If a connected graph contains a single block, we call the graph itself a \emph{block}.

For $\ell \ge 2$ and a finite sequence $G_{1}, \dots, G_{\ell}$ of vertex disjoint graphs, we let the \emph{join} $G_{1} \vee  \cdots \vee G_{\ell}$ be the graph obtained from the disjoint union of $G_{1}, \dots, G_{\ell}$ by joining each vertex in $G_{i}$ to all vertices in $G_{i + 1}$ for $i \in[\ell - 1]$. If $V(G_{i}) = \{x\}$, then we write $G_{1} \vee \cdots \vee G_{i - 1} \vee x \vee G_{i + 1} \vee \cdots \vee G_{\ell}$. Moreover, for vertex disjoint graphs $G_1$ and $G_2$ and for a subgraph $H$ of $G_{2}$, the \emph{join} $G_{1} \vee \leftidx{_H}{G}{_{2}}$ is the graph obtained from the disjoint union of $G_{1}$ and $G_{2}$ by joining each vertex in $G_{1}$ to each vertex in $H$.

\subsection{Domination Critical Graphs}
\label{S:critical}

A study of properties of domination critical graphs was initiated by Sumner and Blitch in their classical 1983 paper~\cite{SuBl83}. Among other results, they showed that every connected $3$-$\gamma$-critical graph of even order contains a perfect matching. Wojcicka~\cite{Wo90} subsequently studied hamiltonian properties of domination critical graphs and showed every connected $3$-$\gamma$-critical graph on at least seven vertices is traceable. Favaron et al.~\cite{FaTiZh97}, Flandrin et al.~\cite{FlTiWeZh99} and Tian et al.~\cite{TiWeZh99} proved further that all connected $3$-$\gamma$-critical graphs with minimum degree at least~$2$ are hamiltonian. Motivated in part by these results, Sumner and Wojcicka (Chapter~16 in~\cite{HHS1}) conjectured in 1998 that all $(k - 1)$-connected $k$-$\gamma$-critical graphs are hamiltonian for all $k \ge 4$. However, their conjecture was disproved seven years later by Yuansheng et al.~\cite{YuChXiYoXi05} who constructed a $3$-connected $4$-$\gamma$-critical non-hamiltonian graph containing 13 vertices. On the positive side,  Kaemawichanurat and Caccetta~\cite{KaCa19} proved the Sumner-Wojcicka Conjecture is true if $k = 4$ and the graphs are claw-free.

\subsection{Connected Domination Critical Graphs}
\label{S:Ccritical}

Kaemawichanurat~\cite{thesis} initiated a study of  connected domination critical graphs. Hamiltonian properties of connected domination critical graphs were subsequently studied by Kaemawichanurat, Caccetta and Ananchuen~\cite{KaCaAn18} who showed that every $2$-connected $k$-$\gamma_{c}$-critical graph is hamiltonian for all $k \in [3]$. Further, they constructed $k$-$\gamma_{c}$-critical graphs that are non-hamiltonian for all $k \ge 4$. Recently, Kaemawichanurat and Caccetta~\cite{KaCa19} proved that every $2$-connected $4$-$\gamma_{c}$-critical claw-free graph is hamiltonian, and they constructed $2$-connected $k$-$\gamma_{c}$-critical claw-free graphs that are non-hamiltonian for all $k \ge 5$. For $5 \le k \le 6$, they proved that every $3$-connected $k$-$\gamma_{c}$-critical claw-free graph is hamiltonian.
Recall that $\zeta(G)$ denotes the number of cut-vertices of $G$, and that if the graph $G$ is clear from the context, we simply write $\zeta$ instead of $\zeta(G)$. Kaemawichanurat and Ananchuen~\cite{KaAn19} showed that a connected domination critical graph cannot have too many cut-vertices.

\begin{thm}{\rm (\cite{KaAn19})}
\label{thm mpm}
For $k \ge 2$, every $k$-$\gamma_{c}$-critical graph has at most $k - 2$ cut-vertices, that is, $\zeta \le k - 2$.
\end{thm}

\section{Main Result}
\label{S:main}

Our aim in this paper is to determine a connection between the traceability of a $k$-$\gamma_{c}$-critical graph and the number of cut-vertices in the graph. More precisely, we shall prove the following result.

\begin{thm}
\label{t:main1}
For $k \ge 4$ and $0 \le \zeta \le k - 2$, every $k$-$\gamma_{c}$-critical graph with $\zeta$ cut-vertices has a hamiltonian path if and only if $k - 3 \le \zeta \le k - 2$.
\end{thm}

\section{Preliminary Results}

In this section, we present some preliminary results that we will need to prove our main theorem, namely Theorem~\ref{t:main1}.
The following result is a simple exercise in most graph theory textbooks.

\begin{ob}\label{o:toughoftrace}
Let $G$ be a graph and let $S$ be a nonempty proper subset of $V(G)$. If $G$ is traceable, then $\omega(G - S) \le |S| + 1$.
\end{ob}

By Observation~\ref{o:toughoftrace}, if $S$ is a vertex cut set of a graph $G$ satisfying $|S| +1 < \omega(G - S)$, then $G$ is non-traceable. Kaemawichanurat, Caccetta and Ananchuen~\cite{KaCaAn18} showed that connected domination critical graphs with small connected domination number are hamiltonian.

\begin{thm}{\rm (\cite{KaCaAn18})}
\label{thm hamiltonian}
Every $k$-$\gamma_{c}$-critical graph is hamiltonian for all $k \in [3]$.
\end{thm}

Chen, Sun, and Ma~\cite{ChSuMa04} characterized all $k$-$\gamma_{c}$-critical graphs for $k \in [2]$.

\begin{thm}{\rm (\cite{ChSuMa04})}
\label{thm chen}
A graph $G$ is $1$-$\gamma_{c}$-critical if and only if $G$ is a complete graph. Moreover, a graph $G$ is $2$-$\gamma_{c}$-critical if and only if $\overline{G} = \cup^{k}_{i = 1}K_{1,n_{i}}$ where $k \ge 2$ and $n_{i} \ge 1$ for all $i \in [k]$.
\end{thm}

Chen et al.~\cite{ChSuMa04} also established fundamental properties of $k$-$\gamma_{c}$-critical graphs for $k \ge 2$.
\vskip 5 pt

\begin{lem}{\rm (\cite{ChSuMa04})}
\label{lem 1}
Let $G$ be a $k$-$\gamma_{c}$-critical graph, and let $x$ and $y$ be a pair of non-adjacent vertices of $G$. If $D_{xy}$ is a $\gamma_{c}$-set of $G + xy$, then the following holds. \\ [-26pt]
\begin{enumerate}
\item $k - 2 \le |D_{xy}| \le k - 1$.
\item $D_{xy} \cap \{x, y\} \ne \emptyset$.
\item If $\{x\} = \{x, y\} \cap D_{xy}$, then $N_{G}(y) \cap D_{xy} = \emptyset$.
\end{enumerate}
\end{lem}

Ananchuen~\cite{An07} established the following properties and structural results of $k$-$\gamma_{c}$-critical graphs that possess cut-vertices.

\begin{lem}{\rm (\cite{An07})}
\label{lem 2}
For $k \ge 3$, if $G$ is a $k$-$\gamma_{c}$-critical graph with a cut-vertex $c$ and if $D$ is a CD-set of $G$, then the following holds. \\ [-26pt]
\begin{enumerate}
\item $G - c$ contains exactly two components.
\item If $C_{1}$ and $C_{2}$ are the components of $G - c$, then $G[N_{C_{1}}(c)]$ and $G[N_{C_{2}}(c)]$ are complete.
\item $c \in D$.
\end{enumerate}
\end{lem}

As remarked earlier, Kaemawichanurat and Ananchuen~\cite{KaAn19} showed in Theorem~\ref{thm mpm} that for $k \ge 2$, every $k$-$\gamma_{c}$-critical graph has at most $k - 2$ cut-vertices, that is, $\zeta \le k - 2$. Further, they also characterized the $k$-$\gamma_{c}$-critical graph with exactly $k - 2$ cut-vertices. To state their results, let $\cS$ be a set of stars $G_1, G_2, \ldots, G_{|\cS|}$ where $|\cS| \ge 2$, $G_i \cong K_{1, n_{i}}$ and $V(G_i) = \{s^{i}_{0}, s^{i}_{1}, \ldots, s^{i}_{n_{i}}\}$ where $s^{i}_{0}$ is the center of the star $G_i$ for $i \in [ \, |\cS| \, ]$.  Let
\[
S = \bigcup^{|\cS |}_{i = 1}\{s^{i}_{0}\} \hspace*{0.5cm} \mbox{and} \hspace*{0.5cm}
S' = \bigcup^{|\cS |}_{i = 1}\{s^{i}_{1}, s^{i}_{2}, ..., s^{i}_{n_{i}}\}.
\]
Moreover, let $S''$ be a (possibly empty) set of  isolated vertices. We note that $|S| = |\cS | \ge 2$. Let $T$ be the vertex disjoint union of these stars $G_1, G_2, \ldots, G_{|\cS|}$. Thus, the complement $\barT$ of $T$ is a complete graph obtained by removing the edges from the stars in $\cS$. We are now in a position to describe the following classes of graphs.

\noindent \textbf{The class $\cB_{1}$}. A graph $G$ in the class $\cB_{1}$ is constructed from the complement $\barT$ of $T$ by adding a new vertex~$b$ and joining it to every vertex of $S'$. The vertex $b$ of $G$ is called the \emph{head} of $G$. A graph in the class $\cB_{1}$ is illustrated in Figure~\ref{f:fig1}.

\vskip -0.25 cm
\setlength{\unitlength}{0.8cm}
\begin{figure}[htb]
\begin{center}
\begin{picture}(13, 7.5)
\put(3, 4){\circle*{0.2}}
\put(7.5, 1.7){\circle*{0.2}}
\put(7.5, 2.7){\circle*{0.2}}
\put(5.5, 4){\oval(1, 4.7)}
\put(7.5, 3.9){\oval(1, 5)}
\put(3, 4){\line(1, 1){2}}
\put(3, 4){\line(1, -1){2}}
\put(7, 3){\line(1, 0){1}}
\put(6.1, 4.9){\line(1, 0){0.8}}
\put(6.1, 4.6){\line(1, 0){0.8}}
\put(6.1, 2.4){\line(1, 0){0.8}}
\put(6.1, 2.1){\line(1, 0){0.8}}
\put(5.5, 6){\circle*{0.2}}
\put(5.5, 5){\circle*{0.2}}
\put(7.5, 6){\circle*{0.2}}
\multiput(7, 6)(-0.7,0){3}{\line(1,0){0.4}}
\multiput(7.5, 6)(-0.7,-0.35){3}{\line(-2,-1){0.5}}
\put(5.45, 5.3){\footnotesize$\vdots$}
\put(7.45, 2){\footnotesize$\vdots$}

\put(5.5, 4){\circle*{0.2}}
\put(5.5, 3){\circle*{0.2}}
\put(7.5, 4){\circle*{0.2}}
\multiput(7, 4)(-0.7,0){3}{\line(1,0){0.4}}
\multiput(7.5, 4)(-0.7,-0.35){3}{\line(-2,-1){0.5}}
\put(5.45, 3.3){\footnotesize$\vdots$}
\put(7.45, 5){\footnotesize$\vdots$}

\put(2.8, 3.5){\footnotesize$b$}
\put(5.35, 1.2){\footnotesize$S'$}
\put(8.3, 4.5){\footnotesize$S$}
\put(8.3, 2.1){\footnotesize$S''$}
\end{picture}
\vskip -0.75 cm
\caption{\label{f:fig1} A graph $G$ in the class $\cB_{1}$}
\end{center}
\end{figure}

\noindent \textbf{The class $\cU(k)$}. Let $B$ be a graph in the class $\cB_{1}$ defined earlier. A graph $G$ in the class $\cU(k)$ is constructed from the graph $B$ and a path $P_{k - 2} \colon c_{0} c_{1} \ldots c_{k - 3}$ of order $k-2$ by joining $c_{k - 3}$ to $b$. A graph $G$ in the class $\cU(k)$ is illustrated by Figure~\ref{f:fig2}.

\vskip -0.25 cm
\setlength{\unitlength}{0.8cm}
\begin{figure}[htb]
\begin{center}
\begin{picture}(13, 5)
\put(0, 3){\circle*{0.2}}
\put(2, 3){\circle*{0.2}}
\put(4, 3){\circle*{0.2}}
\put(8, 3){\circle*{0.2}}
\put(10, 3){\circle*{0.2}}
\put(10.9, 3){\circle{3}}
\put(0, 3){\line(1, 0){10}}

\put(0, 2.5){\footnotesize$c_{0}$}
\put(2, 2.5){\footnotesize$c_{1}$}
\put(4, 2.5){\footnotesize$c_{2}$}
\put(6, 2.5){\footnotesize$\ldots$}
\put(8, 2.5){\footnotesize$c_{k - 3}$}
\put(9.8, 2.5){\footnotesize$b$}

\put(10.3, 4.1){\footnotesize$B \in \cB_{1}$}

\end{picture}
\vskip -1.7 cm
\caption{\label{f:fig2} A graph $G$ in the class $\cU(k)$}
\end{center}
\end{figure}

\medskip
We are now in a position to state the characterization of $k$-$\gamma_{c}$-critical graphs with $k - 2$ cut-vertices.

\begin{thm}{\rm (\cite{KaAn19})}
\label{thm mpm}
For $k \ge 2$, if $G$ is a $k$-$\gamma_{c}$-critical graph, then $\zeta \le k - 2$. Moreover, $\zeta = k - 2$ if and only if $G \in \cU(k)$.
\end{thm}

In order to present the characterization due to Kaemawichanurat~\cite{Ka19} of $k$-$\gamma_{c}$-critical graphs with $\zeta = k - 3$ cut-vertices, we describe next some additional classes of graphs. Let $\textbf{\emph{i}} = (i_{1}, i_{2}, \ldots, i_{k - 3})$ be a $(k - 3)$-tuple such that $i_{1}, i_{2}, \ldots, i_{k - 3} \in \{0, 1\}$ and
$\sum^{k - 3}_{j = 1}i_{j} = 1$.
Thus, there is exactly one $\ell \in [k-3]$ such that $i_{\ell} = 1$ and $i_{\ell'} = 0$ for all $\ell' \in [k-3] \setminus \{\ell\}$.

\noindent \textbf{The class $\cG_{1}(i_{1}, i_{2}, \ldots, i_{k - 3})$}. For a $(k - 3)$-tuple $\textbf{\emph{i}} = (0, 0, \ldots, i_{\ell}, \ldots, 0)$ where $i_{\ell} = 1$ and $i_{\ell'} = 0$ for $1 \le \ell \le k - 4$ and $1 \le \ell' \le k - 3$ where $\ell \ne \ell'$, a graph $G$ in the class $\cG_{1}\textbf{\emph{i}}$ can be constructed from the vertex disjoint paths $c_{0} c_{1} \ldots c_{\ell - 1}$ and $c_{\ell} c_{\ell + 1} \ldots c_{k - 4}$, a copy of a complete graph $K_{n_{\ell}}$ and a block $B \in \cB_{1}$ by adding edges according the join operations
\[
c_{\ell - 1} \vee K_{n_{\ell}} \vee c_{\ell} \hspace*{0.5cm} \mbox{and} \hspace*{0.5cm}
c_{k - 4} \vee b
\]
where $b$ is the head of $B$. Thus, the vertices $c_{\ell - 1}$ and $c_{\ell}$ are joined to every vertex in the complete graph $K_{n_{\ell}}$, and the vertices $c_{k - 4}$ and $b$ are joined. Two examples of graphs in this case when $1 \le \ell \le k - 4$ are illustrated by Figure~\ref{f:fig3} and Figure~\ref{f:fig4}.

\vskip -1.25 cm
\setlength{\unitlength}{0.8cm}
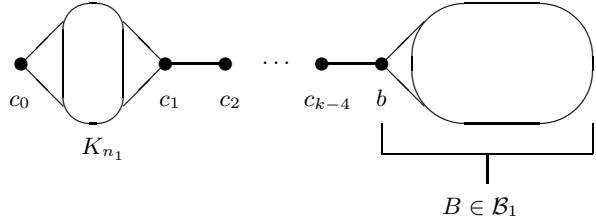
\begin{figure}[htb]
\begin{center}
\begin{picture}(15, 6)
\put(3, 3){\circle*{0.2}}
\put(4.2, 3){\oval(1, 2)}
\put(11, 3){\oval(3, 2)}

\put(3, 3){\line(1, 1){0.7}}
\put(3, 3){\line(1, -1){0.7}}
\put(5.4, 3){\line(-1, 1){0.7}}
\put(5.4, 3){\line(-1, -1){0.7}}

\put(5.4, 3){\line(1, 0){1}}
\put(8, 3){\line(1, 0){1}}

\put(9, 3){\line(1, 1){0.7}}
\put(9, 3){\line(1, -1){0.7}}

\put(9, 1.5){\line(1, 0){3.5}}
\put(9, 1.5){\line(0, 1){0.5}}
\put(12.5, 1.5){\line(0, 1){0.5}}
\put(10.75, 1.5){\line(0, -1){0.5}}

\put(5.4, 3){\circle*{0.2}}
\put(6.4, 3){\circle*{0.2}}
\put(8, 3){\circle*{0.2}}
\put(9, 3){\circle*{0.2}}
\put(7, 3){\footnotesize$\ldots$}
\put(2.8, 2.3){\footnotesize$c_{0}$}
\put(5.3, 2.3){\footnotesize$c_{1}$}
\put(6.3, 2.3){\footnotesize$c_{2}$}
\put(7.7, 2.3){\footnotesize$c_{k - 4}$}
\put(8.9, 2.3){\footnotesize$b$}

\put(4, 1.5){\footnotesize$K_{n_{1}}$}
\put(10, 0.5){\footnotesize$B \in \cB_{1}$}

\end{picture}
\vskip -0.5 cm
\caption{\label{f:fig3}
A graph $G$ in the class $\cG_{1}(i_{1} = 1, 0, 0, \ldots, 0)$}
\end{center}
\end{figure}

\vskip -0.5 cm
\setlength{\unitlength}{0.8cm}
\begin{figure}[htb]
\begin{center}
\begin{picture}(12, 5)
\put(2.5, 3){\circle*{0.2}}
\put(1.5, 3){\circle*{0.2}}
\put(0.5, 3){\circle*{0.2}}

\put(5.2, 3){\oval(1, 2)}
\put(11, 3){\oval(3, 2)}

\put(0.5, 3){\line(1, 0){2}}
\put(4, 3){\line(1, 1){0.7}}
\put(4, 3){\line(1, -1){0.7}}
\put(6.4, 3){\line(-1, 1){0.7}}
\put(6.4, 3){\line(-1, -1){0.7}}

\put(8, 3){\line(1, 0){1}}
\put(9, 3){\line(1, 1){0.7}}
\put(9, 3){\line(1, -1){0.7}}

\put(9, 1.5){\line(1, 0){3.5}}
\put(9, 1.5){\line(0, 1){0.5}}
\put(12.5, 1.5){\line(0, 1){0.5}}
\put(10.75, 1.5){\line(0, -1){0.5}}

\put(4, 3){\circle*{0.2}}
\put(6.4, 3){\circle*{0.2}}
\put(8, 3){\circle*{0.2}}
\put(9, 3){\circle*{0.2}}
\put(7, 3){\footnotesize$\ldots$}
\put(3.1, 3){\footnotesize$\ldots$}

\put(0.3, 2.3){\footnotesize$c_{0}$}
\put(1.3, 2.3){\footnotesize$c_{1}$}
\put(2.3, 2.3){\footnotesize$c_{2}$}

\put(3.8, 2.3){\footnotesize$c_{\ell - 1}$}
\put(6.3, 2.3){\footnotesize$c_{\ell}$}
\put(7.7, 2.3){\footnotesize$c_{k - 4}$}
\put(8.9, 2.3){\footnotesize$b$}

\put(5, 1.5){\footnotesize$K_{n_{\ell}}$}
\put(10, 0.5){\footnotesize$B \in \cB_{1}$}
\end{picture}
\vskip -0.5 cm
\caption{\label{f:fig4}
A graph $G$ in the class $\cG_{1}(0, 0, \ldots, i_{\ell} = 1, 0, \ldots, 0)$}
\end{center}
\end{figure}
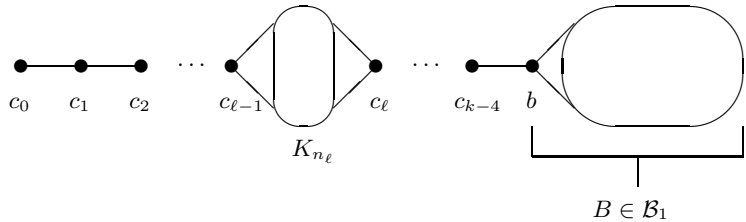

\bigskip
Further, for a $(k - 3)$-tuple $\textbf{\emph{i}} = (0, 0, \ldots, 1)$ where $i_{k - 3} = 1$ and $i_{\ell'} = 0$ for $\ell' \in [k-4]$, a graph $G$ in the class $\cG_{1}\textbf{\emph{i}}$ can be constructed from a path $c_{0} c_{1} \ldots c_{k - 4}$, a copy of a complete graph $K_{n_{k - 3}}$ and a block $B \in \cB_{1}$ by adding edges according the join operation $c_{k - 4} \vee K_{n_{k - 3}} \vee b$, where $b$ is the head of $B$. Thus, the vertices $c_{k-4}$ and $b$ are joined to every vertex in the complete graph $K_{n_{k - 3}}$. An example of a graph in this case is illustrated in Figure~\ref{f:fig5}.

\vskip -0.5 cm
\setlength{\unitlength}{0.8cm}
\begin{figure}[htb]
\begin{center}
\begin{picture}(11, 5)
\put(2.5, 3){\circle*{0.2}}
\put(1.5, 3){\circle*{0.2}}
\put(0.5, 3){\circle*{0.2}}

\put(5.2, 3){\oval(1, 2)}
\put(8.3, 3){\oval(3, 2)}

\put(0.5, 3){\line(1, 0){2}}
\put(4, 3){\line(1, 1){0.7}}
\put(4, 3){\line(1, -1){0.7}}
\put(6.4, 3){\line(-1, 1){0.7}}
\put(6.4, 3){\line(-1, -1){0.7}}

\put(6.3, 3){\line(1, 1){0.7}}
\put(6.3, 3){\line(1, -1){0.7}}

\put(6.3, 1.5){\line(1, 0){3.5}}
\put(6.3, 1.5){\line(0, 1){0.5}}
\put(9.8, 1.5){\line(0, 1){0.5}}
\put(8.05, 1.5){\line(0, -1){0.5}}

\put(4, 3){\circle*{0.2}}
\put(6.4, 3){\circle*{0.2}}
\put(3.1, 3){\footnotesize$\ldots$}

\put(0.3, 2.3){\footnotesize$c_{0}$}
\put(1.3, 2.3){\footnotesize$c_{1}$}
\put(2.3, 2.3){\footnotesize$c_{2}$}

\put(3.8, 2.3){\footnotesize$c_{k - 4}$}
\put(6.3, 2.3){\footnotesize$b$}

\put(5, 1.5){\footnotesize$K_{n_{k - 3}}$}
\put(7.3, 0.5){\footnotesize$B \in \cB_{1}$}
\end{picture}
\vskip -0.5 cm
\caption{\label{f:fig5}
A graph $G$ in the class $\cG_{1}(0, 0, \ldots, 1)$}
\end{center}
\end{figure}
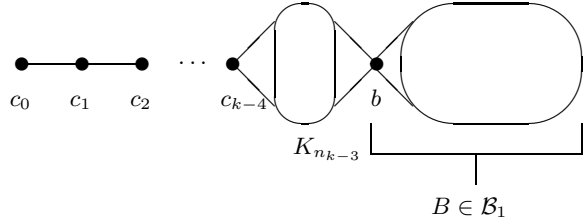

We proceed further by defining a special class of end blocks.

\noindent \textbf{The class $\cB_{2}$}. Let $H$ be a block graph, and so $H$ is a connected graph that contains a single block. The block $H$ belongs to the family $\cB_{2}$ if $\gamma_c(H) = 3$ and $H$ has the following properties. \\ [-26pt]
\begin{enumerate}
\item The block $H$ contains a vertex $b$ such that $N_{H}(b)$ is a complete graph. 
\item Every vertex $v$ of $H$ different from $b$ belongs to some $\gamma_{c}$-set of $H$ of size~$3$.
\item For every pair of non-adjacent vertices $x$ and $y$ in $H - b$, there exists a $\gamma_{c}$-set of $H+xy$ of size~$2$ that contains a neighbor of $b$ in $H$ and contains at least one of $x$ and $y$.
\end{enumerate}
\vskip -0.25 cm
The vertex $b$ is called the \emph{head} of the block $H \in \cB_2$. We note that in property~(b) defined above, the  $\gamma_{c}$-set of $H$ that contains the vertex~$v \in V(H) \setminus \{b\}$ must contain a neighbor of $b$ in $H$ in order to dominate the vertex~$b$.

\noindent \textbf{The class $\cG_{2}(k)$ for $k \ge 5$}. A graph $G$ belongs to the class $\cG_{2}(k)$ for $k \ge 5$ if it can be constructed from the vertex disjoint union of a path $c_{0} c_{1} \ldots c_{k - 4}$ and a block graph $H \in \cG_{2}$ with head $b$ by adding the edge $bc_{k - 4}$.

We are now in a position to state the characterization of $k$-$\gamma_{c}$-critical graphs with $\zeta = k - 3$ cut-vertices due to Kaemawichanurat~\cite{Ka19}.

\begin{thm}{\rm (\cite{Ka19})}
\label{thm k-3}
For $k \ge 4$, if $G$ is a $k$-$\gamma_{c}$-critical graph with $k - 3$ cut-vertices, then $G \in \cG_{1}(i_{1}, i_{2}, \ldots, i_{k - 3}) \cup \cG_{2}(k)$.
\end{thm}

\section{Traceability of $k$-$\gamma_{c}$-Critical Graphs}
\label{S:traceability}

In this section, we show that, for $k \ge 4$ and $k - 3 \le \zeta \le k - 2$, every $k$-$\gamma_{c}$-critical graph with $\zeta$ cut-vertices contains a hamiltonian path. We first prove basic properties of $k$-$\gamma_{c}$-critical graphs.

In what follows, let $B$ be a graph in the class $\cB_{1}$ of order $n_{0}$ and let the vertex $b$ be the head of $B$. For notational convenience, we sometimes rename the vertex $b$ as the vertex $c_{k - 2}$. We show first that there exists a hamiltonian path in $B$ that contains the vertex $b$ as one of its ends.

\begin{lem}
\label{l:Hampath}
If $B \in \cB_{1}$ with the vertex $b$ as its head, then there exists a hamiltonian path $P_{B}$ of $B$ having $b$ as one of its ends.
\end{lem}
\proof  By the construction of the graph $B \in \cB_{1}$, we have $S' = N_{B}(b)$ and $S \cup S'' = V(B) \setminus N_{B}[b]$. Further, we note that $B[S' \cup \{b\}]$ and $B[S \cup S'']$ are complete subgraphs. Since $|S| \ge 2$, every vertex in $S'$ has at least one neighbor in $S$. Let $uv$ be an arbitrary edge in $G$ where $u \in S'$ and $v \in S$. Further, let $u'$ be an arbitrary vertex in $S'$ different from~$u$. Let $P_u$ be a hamiltonian path in $G[S']$ that starts at the vertex $u'$ and ends at the vertex~$u$. Let $P_v$ be a hamiltonian path in $G[S \cup S'']$ that starts at the vertex~$v$. Let $P_B$ be the hamiltonian path of $B$ that starts at the vertex $b$, proceeds along the edge $bu'$ to $u'$, follows the hamiltonian path $P_u$ from $u'$ to $u$, proceeds along the edge $uv$ to $v$, and then follows the hamiltonian path $P_v$ starting at the vertex~$v$. By construction, the hamiltonian path $P_B$ of $B$ has the vertex $v$ as one of its ends.~\qed

By Lemma~\ref{l:Hampath}, there exists a hamiltonian path $P_{B}$ of $B$ having $b$ as one of its ends. Let $b'$ be the other end of the path $P_B$. As a consequence of Theorem \ref{thm mpm} and Lemma~\ref{l:Hampath}, we obtain the following lemma.

\begin{lem}\label{lem Tk-2}
If $G$ is a $k$-$\gamma_{c}$-critical graph with $k - 2$ cut-vertices, then $G$ is traceable.
\end{lem}
\proof Let $G$ be a $k$-$\gamma_{c}$-critical graph with $k - 2$ cut-vertices. By Theorem~\ref{thm mpm}, the graph $G \in \cU(k)$. Therefore, $G$ is constructed from a graph $B \in \cB_{1}$ with head $b$ and a path $P_{k - 2} \colon c_{0} c_{1} \ldots c_{k - 3}$ by joining $c_{k - 3}$ to $b$. The path $c_{0} c_{1} \ldots c_{k - 3}$ can therefore be extended to a hamiltonian path of $G$ by proceeding along the edge $c_{k - 3}b$ from $c_{k - 3}$ to $b$, and then following the hamiltonian path $P_{B}$ from $b$ to $b'$ to yield the hamiltonian path $c_{0}  c_{1}  \ldots c_{k - 3}  b P_{B} b'$ of $G$.~\qed

We show next that every graph in the class $\cG_{1}(i_{1}, i_{2}, \ldots, i_{k - 3})$ has a hamiltonian path.

\begin{lem}\label{lem Tk-3G1}
If $G \in \cG_{1}(i_{1}, i_{2}, \ldots, i_{k - 3})$, then $G$ is traceable.
\end{lem}
\proof Suppose that $G \in \cG_{1}(i_{1}, i_{2}, \ldots, i_{k - 3})$, where $i_{\ell} = 1$ and $i_{\ell'} = 0$ for $\ell, \ell' \in [k - 4]$ and $\ell \ne \ell'$. Let $Q$ be a hamiltonian path in the copy of $K_{n_{\ell}}$ used in the construction of $G$, and let $w_1$ and $w_{n_{\ell}}$ be the start and final vertex of the path $Q$.

We first consider the case when $\ell = 1$, and so $G \in \cG_{1}(1, 0, 0, \ldots, 0)$. The path that starts at the vertex $c_0$, proceeds along the edge $c_0w_1$ to $w_1$, follows the hamiltonian path $Q$ from $w_1$ to $w_{n_{1}}$, proceeds along the edge $w_{n_{1}}c_1$ to $c_1$, follows the path $c_1 c_2 \ldots c_{k-4}$, proceeds along the edge $c_{k-4}b$ to $b$, and then follows the hamiltonian path $P_{B}$ from $b$ to $b'$ yield the hamiltonian path $c_{0} P' c_{1}  \ldots c_{k - 4}  b P_{B} b'$ of $G$.

Secondly we consider the case when $2 \le \ell \le k - 3$. We note that $c_{\ell - 1} \succ K_{n_{\ell}}$ and $c_{\ell} \succ K_{n_{\ell}}$. Starting with the path $c_0 c_1 \ldots c_{\ell - 1}$ from $c_0$ to $c_{\ell - 1}$, we proceed along the edge $c_{\ell - 1}w_1$ from $c_{\ell - 1}$ to $w_1$, follow the hamiltonian path $Q$ from $w_1$ to $w_{n_{\ell}}$, proceed along the edge $w_{n_{\ell}}c_\ell$ from $w_{n_{\ell}}$ to $c_\ell$, follow the path $c_\ell \ldots c_{k-4} b$ from $c_\ell$ to $b$, and follow the hamiltonian path $P_{B}$ from $b$ to $b'$ to yield the hamiltonian path $c_{0} c_1 \ldots c_{\ell - 1} Q c_{\ell}  \ldots b P_{B} b'$ of~$G$.

Thirdly we consider the case when $\ell = k - 3$, and so $G \in \cG_{1}(0, 0, \ldots, 1)$. In this case, $c_{k-4} \succ K_{n_{k-3}}$ and $b \succ K_{n_{k-3}}$. Starting with the path $c_0 c_1 \ldots c_{k-4}$ from $c_0$ to $c_{k-4}$, we proceed along the edge $c_{k-4}w_1$ from $c_{k-4}$ to $w_1$, follow the hamiltonian path $Q$ from $w_1$ to $w_{n_{k-3}}$, proceed along the edge $w_{n_{k-3}}b$ from $w_{n_{k-3}}$ to $b$, and follow the hamiltonian path $P_{B}$ from $b$ to $b'$ to yield the hamiltonian path $c_{0} c_1 \ldots c_{k-4} Q b P_{B} b'$ of~$G$. This completes the proof of Lemma~\ref{lem Tk-3G1}.~\qed

\medskip
We are now in a position to prove that all $k$-$\gamma_{c}$-critical graphs with $\zeta$ cut-vertices are traceable when $\zeta \in \{k - 3, k - 2\}$.

\begin{thm}\label{thm traceability}
For $k \ge 4$ and $\zeta \in \{k-3,k-2\}$, if $G$ is a $k$-$\gamma_{c}$-critical graph with $\zeta$ cut-vertices, then $G$ is traceable.
\end{thm}
\proof For $k \ge 4$ and $\zeta \in \{k-3,k-2\}$, let $G$ be a $k$-$\gamma_{c}$-critical graph with $\zeta$ cut-vertices. If $\zeta = k - 2$, then by Lemma~\ref{lem Tk-2}, the graph $G$ is traceable. Hence we may assume that $\zeta = k - 3$, for otherwise the desired result follows. By Theorem~\ref{thm k-3}, $G \in \cG_{1}(i_{1}, i_{2}, \ldots, i_{k - 3}) \cup \cG_{2}(k)$. If $G \in \cG_{1}(i_{1}, i_{2}, \ldots, i_{k - 3})$, then, by Lemma~\ref{lem Tk-3G1}, the graph $G$ is traceable. Hence we may assume that $G \in \cG_{2}(k)$, for otherwise the desired result follows. Thus, $k \ge 5$ and $G$ can be constructed from the vertex disjoint union of a path $P \colon c_{0} c_{1} \ldots c_{k - 4}$ and a block graph $H \in \cG_{2}$ with head $b$ by adding the edge $bc_{k - 4}$. Let
\[
A = N_{H}(b) \hspace*{0.5cm} \mbox{and} \hspace*{0.5cm} \barA  = V(H) \setminus N_{H}[b].
\]
By construction of the graph $H \in \cG_{2}$, we note that $G[A]$ is a complete subgraph. We now consider $G[\, \barA \, ]$.  Let $P^{1} \colon x^{1}_{1} x^{1}_{2} \ldots x^{1}_{n_{1}}$ be a longest path in $G[\, \barA \, ]$. We note that $P^{1}$ is a subgraph of $G[\, \barA \, ]$ and thus, $x^{1}_{j}$ and $x^{1}_{j'}$ may be adjacent for $1 \le j \le j' + 2 \le n_{1}$. If $\barA_1 = V(P^{1})$ and $\barA \setminus \barA_1 \ne \emptyset$, then we let $P^{2} \colon x^{2}_{1} x^{2}_{2} \ldots x^{2}_{n_{2}}$ be a longest path in $G[\, \barA \setminus \barA_1]$. Continuing in this way, for $i \ge 1$ if the paths $P^1, P^2, \ldots, P^{i}$ are defined and $\barA \setminus \barA_i \ne \emptyset$ where
\[
\barA_i = \bigcup_{j=1}^i V(P^{j}),
\]
then we let $P^{i+1} \colon x^{i+1}_{1} x^{i+1}_{2} \ldots x^{i+1}_{n_{1}}$ be a longest path in $G[\, \barA \setminus \barA_i]$. Continuing in this way, let $z \ge 1$ be the smallest integer such that  $\barA \setminus \barA_z = \emptyset$. Thus either $z = 1$, in which case $\barA = V(P^1)$, or $z \ge 2$, in which case $(V(P^1), V(P^2), \ldots, V(P^z))$ is a partition of $\barA$ where each set $V(P^i)$ is nonempty for all $i \in [z]$. By definition of the paths $P^i$ for $i \in [z]$, we note that
\[
|V(P^{1})| \ge |V(P^{2})| \ge \cdots \ge |V(P^{z})|.
\]

The structure of $G[\, \barA \, ]$ is illustrated in Figure~\ref{f:fig6}.

\setlength{\unitlength}{0.8cm}
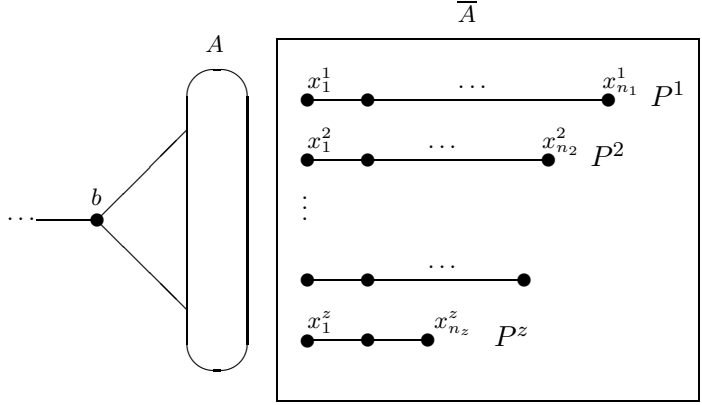
\begin{figure}[htb]
\begin{center}
\begin{picture}(13, 6)
\put(2, 2.5){\circle*{0.2}}
\put(2, 2.5){\line(1, 1){1.5}}
\put(2, 2.5){\line(1, -1){1.5}}
\put(4, 2.5){\oval(1, 5)}
\put(5, -0.5){\framebox(7,6){}}
\put(2, 2.5){\line(-1, 0){1}}
\put(5.5, 4.5){\circle*{0.2}}
\put(5.5, 3.5){\circle*{0.2}}
\put(5.5, 1.5){\circle*{0.2}}
\put(5.5, 0.5){\circle*{0.2}}
\put(6.5, 4.5){\circle*{0.2}}
\put(6.5, 3.5){\circle*{0.2}}
\put(6.5, 1.5){\circle*{0.2}}
\put(6.5, 0.5){\circle*{0.2}}
\put(10.5, 4.5){\circle*{0.2}}
\put(9.5, 3.5){\circle*{0.2}}
\put(9.1, 1.5){\circle*{0.2}}
\put(7.5, 0.5){\circle*{0.2}}
\put(8, 4.7){\footnotesize$\ldots$}
\put(7.5, 3.7){\footnotesize$\ldots$}
\put(7.5, 1.7){\footnotesize$\ldots$}
\put(5.5, 4.5){\line(1, 0){5}}
\put(5.5, 3.5){\line(1, 0){4}}
\put(5.5, 1.5){\line(1, 0){3.6}}
\put(5.5, 0.5){\line(1, 0){2}}
\put(5.5, 4.7){\footnotesize$x^{1}_{1}$}
\put(10.4, 4.7){\footnotesize$x^{1}_{n_{1}}$}
\put(5.5, 3.7){\footnotesize$x^{2}_{1}$}
\put(9.4, 3.7){\footnotesize$x^{2}_{n_{2}}$}
\put(5.5, 0.7){\footnotesize$x^{z}_{1}$}
\put(7.6, 0.7){\footnotesize$x^{z}_{n_{z}}$}
\put(5.4, 2.5){\footnotesize$\vdots$}
\put(0.5, 2.5){\footnotesize$\ldots$}
\put(3.8, 5.3){\footnotesize$A$}
\put(8, 5.8){\footnotesize$\barA $}
\put(11.2, 4.4){$P^{1}$}
\put(10.2, 3.4){$P^{2}$}
\put(8.6, 0.4){$P^{z}$}
\put(1.9, 2.75){\footnotesize$b$}
\end{picture}
\vskip 0.25 cm
\caption{\label{f:fig6}
The structure of $G[ \, \barA \, ]$ in the proof of Theorem~\ref{thm traceability}}
\end{center}
\end{figure}

We proceed further with the following series of claims.

\begin{claim}
\label{lem keylem}
The set $\{x^{1}_{j_{1}}, x^{2}_{j_{2}}, \ldots, x^{z}_{j_{z}}\}$ is an independent set for all $j_{i} \in \{1, n_{i}\}$ and $i \in [z]$.
\end{claim}
\proof Suppose, to the contrary, that $x^{i}_{j_{i}}x^{i'}_{j_{i'}} \in E(G)$ for some $i$ and $i'$ where $1 \le i < i' \le z$ where $j_{i} \in \{1, n_{i}\}$ and $j_{i'} \in \{1, n_{i'}\}$. Renaming the vertices on the path $P^i$ and $P^{i'}$ if necessary, we may assume without loss of generality that $j_{i} = n_i$ and $j_{i'} = 1$. We now consider the path $P^*$ obtained from $P^i$ by proceeding along the edge $x^{i}_{n_i}x^{i'}_{1}$ from $x^{i}_{n_i}$ to $x^{i'}_{1}$, and then following the path $P^{i'}$ from $x^{i'}_{1}$ to $x^{i'}_{n_{i'}}$. If $i = 1$, then $P^*$ is a longer path in $G[\, \barA \, ]$ that $P^1$, contradicting the maximality of $P^1$. If $i \ge 2$, then $P^*$ is a longer path in $G[\, \barA \setminus \barA_{i-1}]$ that $P^i$, contradicting the maximality of $P^i$.~\qed

\medskip
In what follows, we adopt the following notation. If $x$ and $y$ are two non-adjacent vertices of $G$, then we let $D_{xy}$ denote a $\gamma_{c}$-set of $G + xy$.

\begin{claim}\label{lem p2}
If $x$ and $y$ are two non-adjacent vertices of $G[\, \barA \,]$, then $|D_{xy} \cap (A \cup \barA)| = 2$, implying that $|D_{xy} \cap \{x, y\}| = 1$ and $|D_{xy} \cap A| = 1$.
\end{claim}
\proof Let $x$, $y$ and $D_{xy}$ be as defined in the statement of the claim. We now consider the graph $G + xy$. Since $G$ is a $k$-$\gamma_{c}$-critical graph, By Lemma~\ref{lem 1}(a) implies that $|D_{xy}| \le k - 1$. Further, Lemma~\ref{lem 1}(b) implies that $D_{xy} \cap \{x, y\} = 1$. Renaming $x$ and $y$ if necessary, we may assume that $x \in D_{xy}$. If $c_{0} \notin D_{xy}$, then $c_{1} \in D_{xy}$ to dominate $c_{0}$. If $c_{0} \in D_{xy}$, then, since the subgraph, $(G + xy)[D_{xy}]$, of $G + xy$ induced by the set $D_{xy}$ is connected and since $c_1$ is the only neighbor of $c_0$ in $G + xy$, we must have $c_{1} \in D_{xy}$. Hence, in both cases, $c_{1} \in D_{xy}$. Recall that $x \in D_{xy} \cap \barA$. Since $(G + xy)[D_{xy}]$ is a connected graph that contains both $c_1$ and $x$, the structure of the graph $G$ implies that $D_{xy}$ contains all vertices of the path $P \colon c_{0} c_{1} \ldots c_{k - 4}$ except possibly for the vertex $c_0$, the vertex $b$, at least one neighbor of $b$ in $A$, and at least one vertex in $\barA$, namely the vertex~$x$. Thus, $D_{xy}$ contains at least~$(|V(P)| - 1) + 3 = (k-4) + 3 = k-1$ vertices, and so $|D_{xy}| \ge k-1$. As observed earlier, by Lemma~\ref{lem 1}(a) we have $|D_{xy}| \le k - 1$. Consequently, $|D_{xy}| = k-1$, implying that $D_{xy} = (V(P) \setminus \{c_0\}) \cup \{b,u,x\}$, where $u \in A$ and $ux \in E(G)$. In particular, we note that $|D_{xy} \cap V(H - b)| = |D_{xy} \cap (A \cup \barA)| = |\{u,x\}| = 2$, and so $|D_{xy} \cap \{x, y\}| = 1$ and $|D_{xy} \cap A| = 1$.~\qed

\medskip
In what follows, for notational convenience we let
$A_{0} = V(P) \cup \{b\} \cup A$, and so $A_0 = V(G) \setminus \barA$.

\begin{claim}\label{lem p3}
If $R$ is a proper subset of vertices of $A$, where possibly $R = \emptyset$, and $v$ is an arbitrary vertex in $A \setminus R$, then there exists a path $P_{R,v}$ from $c_{0}$ to $v$ containing every vertex in $A_{0} \setminus R$.
\end{claim}
\proof Recall that $A = N_B(H)$ and $G[A]$ is a complete graph. Since $R \subset A$, we note therefore that $G[A - R]$ is a complete subgraph. Let $P_v$ be a hamiltonian path in $G[A - R]$ that ends at the vertex~$v$, and let $v'$ be the start vertex of $P_v$ (possibly, $v = v'$). The path $P_{R,v}$ that starts at the vertex $c_0$, follows the path $P$ to $c_{k-4}$, proceeds along the edge $c_{k-4}b$ from $c_{k-4}$ to $b$, along the edge $bv'$ from $b$ to $v'$, and then follows the path $P_v$ is a path from $c_{0}$ to $v$ containing every vertex in $A_{0} \setminus R$.~\qed

\begin{claim}
\label{cl 1}
If $z = 1$, then $G$ is traceable.
\end{claim}
\proof Suppose that $z = 1$, and so $\barA = V(P^1)$. Suppose that $x^{1}_{1}$ or $x^{1}_{n_{1}}$ is adjacent to some vertex $y$ of $A$. Renaming vertices if necessary, we may assume that $y$ is adjacent to $x^{1}_{1}$. By Claim~\ref{lem p3} with $R = \emptyset$, there exists a path $P_{R,y}$ from $c_{0}$ to $y$ containing every vertex in $A_{0}$. The path $P_{R,y}$ can be extended to a hamiltonian path of $G$ by proceeding along the edge $yx^{1}_{1}$ from $y$ to $x^{1}_{1}$, and then following the path $P^1$ from $x^{1}_{1}$ to $x^{1}_{n_{1}}$. Thus, we may assume that neither $x^{1}_{1}$ nor $x^{1}_{n_{1}}$ is adjacent to any vertex of $A$, for otherwise $G$ is traceable as desired. Since $H \in \cB_2$ is a connected graph, this implies that $|V(P^{1})| \ge 3$.

We show next that $x^{1}_{1}x^{1}_{n_{1}} \in E(G)$.
Suppose, to the contrary, that $x^{1}_{1}x^{1}_{n_{1}} \notin E(G)$. In this case, we consider $G + x^{1}_{1}x^{1}_{n_{1}}$. For notational simplicity, let $D^* = D_{x^{1}_{1}x^{1}_{n_{1}}}$. By Claim~\ref{lem p2}, we have $|D^* \cap \barA| = |D^* \cap A| = 1$. Further, $|D^* \cap \{x^{1}_{1}, x^{1}_{n_{1}}\}| = 1$. Renaming $x^{1}_{1}$ and $x^{1}_{n_{1}}$ if necessary, we may assume that $D^* \cap \{x^{1}_{1}, x^{1}_{n_{1}}\} = \{x^{1}_{1}\}$. Let $D^* \cap A = \{u\}$. By the connectedness of $(G + x^{1}_{1}x^{1}_{n_{1}})[D^*]$, this implies that $x^{1}_{1}u \in E(G)$, contradicting our earlier assumption that $x^{1}_{1}$ is not adjacent to any vertex in $A$. Hence, $x^{1}_{1}x^{1}_{n_{1}} \in E(G)$.

Since $x^{1}_{1}x^{1}_{n_{1}} \in E(G)$, we note that $C \colon P^1 + x^{1}_{1}x^{1}_{n_{1}}$ is a hamiltonian cycle of $G[\, \barA \,]$. Since $G$ is a connected graph, there exists a vertex $v$ in $A$ which is adjacent to a vertex of $P^{1}$, say to $x^{1}_{j}$ for some $j$ where $1 < j < n_{1}$. By Claim~\ref{lem p3} with $R = \emptyset$, there exists a path $P_{R,v}$ from $c_{0}$ to $v$ containing every vertex in $A_{0}$. The path $P_{R,y}$ can be extended to a hamiltonian path of $G$ by proceeding along the edge $vx^{1}_{j}$ from $v$ to $x^{1}_{j}$, and then following a hamiltonian path in the cycle $C$ starting at the vertex $x^{1}_{j}$. Thus, $G$ is traceable.~\qed

\newpage
\begin{claim}
\label{cl 2}
If $z = 2$, then $G$ is traceable.
\end{claim}
\proof Suppose that $z = 2$, and so $\barA = V(P^1) \cup V(P^{2})$. Recall that
$|V(P^{1})| \ge |V(P^{2})|$. By Claim~\ref{lem keylem}, the vertex $x^{1}_{1}$ (respectively, $x^{1}_{n_1}$) is adjacent to neither $x^{2}_{1}$ nor  $x^{2}_{n_2}$. In particular, $x^{1}_{1}x^{2}_{1} \notin E(G)$. We now consider the graph $G + x^{1}_{1}x^{2}_{1}$. For notational simplicity, let $D_{1,2} = D_{x^{1}_{1}x^{2}_{1}}$. By Claim~\ref{lem p2}, we have $|D_{1,2} \cap \barA| = |D_{1,2} \cap A| = 1$. Further, $|D_{1,2} \cap \{x^{1}_{1}, x^{2}_{1}\}| = 1$. Let $D_{1,2} \cap \barA = \{y\}$. We consider the cases $x^{1}_{1} \in D_{1,2}$ and $x^{2}_{1} \in D_{1,2}$ separately.

\begin{subclaim}
\label{lem z201}
If $x^{1}_{1} \in D_{1,2}$, then $G$ is traceable.
\end{subclaim}
\proof Suppose that $x^{1}_{1} \in D_{1,2}$. Thus in this case, $D_{1,2} \cap (A \cup \barA) = \{x^{1}_{1},y\}$. Since $(G + x^{1}_{1}x^{2}_{1})[D_{1,2}]$ is a connected graph, $x^{1}_{1}y \in E(G)$. Moreover by Lemma~\ref{lem 1}(c), $yx^{2}_{1} \notin E(G)$.

\begin{unnumbered}{Claim~\ref{lem z201}.1}
If $|V(P^{1})| \le 2$, then $G$ is traceable.
\end{unnumbered}
\proof Suppose that $|V(P^{1})| \le 2$. Suppose that $|V(P^{1})| = 1$. Since $|V(P^{1})| \ge |V(P^{2})|$, we therefore have $|V(P^{2})| = 1$, and so $P^1$ and $P^2$ consists of the single vertices $x^{1}_{1}$ and $x^{2}_{1}$, respectively. But then $\{y_1,y_2\}$ is a CD-set of $H$, where $y_i$ is an arbitrary neighbor of $x^{i}_{1}$ that belongs to $A$ for $i \in [2]$, and so $\gamma_c(H) \le 2$, contradicting the fact that $\gamma_{c}(H) = 3$. Hence, $|V(P^{1})| = 2$, and so $n_1 = 2$ and $P^1$ is the path $x^{1}_{1} x^{1}_{2}$. As observed earlier, $|V(P^{2})| \le 2$.

Suppose firstly that $|V(P^{2})| = 1$. Thus, $P^2$ consists of the single vertex $x^{2}_{1}$, and $\barA  = \{x^{1}_{1}, x^{1}_{2}, x^{2}_{1}\}$. By Claim~\ref{lem keylem}, the vertex $x^{2}_{1}$ is adjacent to neither $x^{1}_{1}$ nor $x^{1}_{2}$. Let $u$ be an arbitrary neighbor of $x^{2}_{1}$ in the connected graph $G$. We note that $u \in A$. If $u = y$, then $\{y, x^{1}_{1}\} \succ_{c} H$, implying that $\gamma_c(H) \le 2$, a contradiction. Thus, $u \ne y$. Since $H$ is a $2$-connected graph, the vertex $x^{1}_{2}$ has a neighbor, $w$ say, different from $x^{1}_{1}$. We note that $w \in A$. If $w \in \{u,y\}$, then $\{u,y\} \succ_{c} H$, a contradiction. Hence, the vertices $u$, $w$ and $y$ are distinct vertices in $A$. By Claim~\ref{lem p3} with $R = \{u,w\}$, there exists a path $P_{R,y}$ from $c_{0}$ to $y$ containing every vertex in $A_{0} \setminus \{u,w\}$. The path $P_{R,y}$ can be extended to a hamiltonian path of $G$ by proceeding along the edge $yx^{1}_{1}$ from $y$ to $x^{1}_{1}$, and then following the path $x^{1}_{1} x^{1}_{2} w u x^{2}_{1}$ from $x^{1}_{1}$ to $x^{1}_{2}$; that is, the path
\[
c_0 P_{R,y} y, x^{1}_{1} x^{1}_{2} w u x^{2}_{1}
\]
is a hamiltonian path in $G$. Hence we may assume that $|V(P^{2})| = 2$, for otherwise $G$ is traceable, as desired. Thus, $\barA  = \{x^{1}_{1}, x^{1}_{2}, x^{2}_{1}, x^{2}_{2}\}$. Recall that $D_{1,2} \cap (A \cup \barA) = \{x^{1}_{1},y\}$ and that the vertex $x^{1}_{1}$ is adjacent to neither $x^{2}_{1}$ nor  $x^{2}_{2}$. Further, $yx^{2}_{1} \notin E(G)$. These observations imply that $yx^{2}_{2} \in E(G)$ in order for $D_{1,2}$ to dominate the vertex $x^{2}_{2}$ in $G + x^{1}_{1}x^{2}_{1}$. By Claim~\ref{lem keylem}, the vertex $x^{1}_{2}$ is adjacent to neither $x^{2}_{1}$ nor $x^{2}_{2}$. Since $H$ is a $2$-connected graph, the vertex $x^{1}_{2}$ has a neighbor, $u$ say, different from $x^{1}_{1}$. We note that $u \in A$. If $u = y$, then $\{y, x^{2}_{2}\} \succ_{c} H$, a contradiction. Hence, $u \ne y$. By Claim~\ref{lem p3} with $R = \{y\}$, there exists a path $P_{R,u}$ from $c_{0}$ to $u$ containing every vertex in $A_{0} \setminus \{y\}$. The path $P_{R,u}$ can be extended to a hamiltonian path of $G$ by proceeding along the edge $ux^{1}_{2}$ from $u$ to $x^{1}_{2}$, and then following the path $x^{1}_{2} x^{1}_{1} y x^{2}_{2} x^{2}_{1}$ from $x^{1}_{2}$ to $x^{2}_{1}$; that is, the path
\[
c_0 P_{R,u} u, x^{1}_{2} x^{1}_{1} y x^{2}_{2} x^{2}_{1}
\]
is a hamiltonian path in $G$.  This completes the proof of Claim~\ref{lem z201}.1.~\smallqed

\medskip
By Claim~\ref{lem z201}.1, we may assume that $|V(P^{1})| \ge 3$, for otherwise $G$ is traceable and the desired result holds.

\begin{unnumbered}{Claim~\ref{lem z201}.2}
If $x^{1}_{1}x^{1}_{n_{1}} \notin E(G)$, then $G$ is traceable.
\end{unnumbered}
\proof Suppose that $x^{1}_{1}x^{1}_{n_{1}} \notin E(G)$. This implies that $|V(P^{1})| \ge 3$. In order to dominate the vertex $x^{1}_{n_{1}}$ in $G + x^{1}_{1}x^{2}_{1}$, we must have that $yx^{1}_{n_{1}} \in E(G)$. We now consider the graph $G + x^{1}_{1}x^{1}_{n_{1}}$. For notational simplicity, let $D_{1,n_1} = D_{x^{1}_{1}x^{1}_{n_1}}$. By Claim~\ref{lem p2}, we have $|D_{1,n_1} \cap \barA| = |D_{1,n_1} \cap A| = 1$. Further, $|D_{1,n_1} \cap \{x^{1}_{1}, x^{1}_{n_1}\}| = 1$. Let $\{u\} = D_{1,n_1} \cap \barA$. As observed earlier, the vertex $x^{2}_{1}$ is adjacent to neither $x^{1}_{1}$ nor $x^{1}_{n_1}$. In order to dominate the vertex $x^{2}_{1}$ in $G + x^{1}_{1}x^{1}_{n_{1}}$, we must have that $ux^{2}_{1} \in E(G)$.

By Claim~\ref{lem p3} with $R = \{u\}$, there exists a path $P_{R,y}$ from $c_{0}$ to $y$ containing every vertex in $A_{0} \setminus \{u\}$. By the connectedness of $(G + x^{1}_{1}x^{1}_{n_{1}})[D_{1,n_1}]$, the vertex $u$ is adjacent to the vertex in $D_{1,n_1} \cap \{x^{1}_{1}, x^{1}_{n_{1}}\}$. Renaming the vertices $x^{1}_{1}$ and $x^{1}_{n_{1}}$ if necessary, we may assume without loss of generality that $x^{1}_{n_{1}} \in D_{1,n_1}$. With this assumption, $ux^{1}_{n_{1}} \in E(G)$. The path $P_{R,y}$ can be extended to a hamiltonian path of $G$ by proceeding along the edge $yx^{1}_{1}$ from $y$ to $x^{1}_{1}$, following the path $P^1$ from $x^{1}_{1}$ to $x^{1}_{n_{1}}$, proceeding along the edge $x^{1}_{n_{1}}u$ from $x^{1}_{n_{1}}$ to $u$, proceeding along the edge $ux^{2}_{1}$ from $u$ to $x^{2}_{1}$, and then following the path $P^2$ from $x^{2}_{1}$ to $x^{2}_{n_{2}}$; that is, the path
\[
c_0 P_{R,y} y, x^{1}_{1} P^1 x^{1}_{n_{1}}, u , x^{2}_{1} P^2 x^{2}_{n_{2}}
\]
is a hamiltonian path in $G$.~\smallqed

\medskip
By Claim~\ref{lem z201}.2, we may assume that $x^{1}_{1}x^{1}_{n_{1}} \in E(G)$, for otherwise $G$ is traceable and the desired result holds. Since $x^{1}_{1}x^{1}_{n_{1}} \in E(G)$, we note that $C \colon P^1 + x^{1}_{1}x^{1}_{n_{1}}$ is a cycle in $G[\, \barA \,]$.

\begin{unnumbered}{Claim~\ref{lem z201}.3}
If $|V(P^{2})| = 1$, then $G$ is traceable.
\end{unnumbered}
\proof Suppose that $|V(P^{2})| = 1$. Thus, $P^2$ consists of the single vertex $x^2_1$. By Claim~\ref{lem keylem}, the vertex $x^{2}_{1}$ is adjacent to neither $x^{1}_{1}$ nor $x^{1}_{n_1}$. If $x^{2}_{1}$ is adjacent to $x^{1}_{j}$ for some $1 < j < n_{1}$, then the $(x^{1}_{j+1},x^{1}_{j})$-path on $C$ that does not contain the edge $x^{1}_{j}x^{1}_{j+1}$ can be extended to a longer path in $G[\, \barA \,]$ by adding to it the vertex $x^{2}_{1}$ and the edge $x^{1}_{j}x^{2}_{1}$, contradicting the maximality of the path $P^1$. Hence, the vertex $x^{2}_{1}$ is adjacent to no vertex of $P^{1}$. By the connectivity of $G$ and the maximality of the path $P^2$, the vertex $x^{2}_{1}$ is adjacent to a vertex, $u$ say, in $A$.

If $y \succ P^{1}$, then $\{y, u\} \succ_{c} H$, implying that $\gamma_c(H) \le 2$, a contradiction. Thus, $y$ does not dominate $P^{1}$. Let $j$ be the smallest integer so that $yx^{1}_{j}$ is not an edge. Since $x^{1}_{1}y \in E(G)$, we note that $j \in [n_1] \setminus \{1\}$. By the choice of $j$, we note that $yx^{1}_{\ell} \in E(G)$ for all $\ell \in [j-1]$.

We now consider the graph $G + x^{1}_{j}x^{2}_{1}$. For notational simplicity, let $D^* = D_{x^{1}_{j}x^{2}_{1}}$. By Claim~\ref{lem p2}, we have $|D^* \cap \barA| = |D^* \cap A| = 1$. Further, $|D^* \cap \{x^{1}_{j},x^{2}_{1}\}| = 1$. Let $\{w\} = D^* \cap A$. Since $y$ is adjacent to neither $x^{1}_{j}$ nor $x^{2}_{1}$, we note that $w \ne y$. If $x^{2}_{1} \in D^*$, then since $x^{2}_{1}$ is not adjacent to any vertex of $P^{1}$, we note that $w \succ P^{1} - x^{1}_{j}$. implying that $\{w, x^{1}_{j - 1}\} \succ_{c} H$, and so $\gamma_c(H) \le 2$, a contradiction. Hence, $x^{2}_{1} \notin D^*$, implying that $x^{1}_{j} \in D^*$; that is, $D^* \cap \barA = \{ x^{1}_{j} \}$. By Lemma~\ref{lem 1}(c), we note that $wx^{2}_{1} \notin E(G)$. Since $(G + x^{1}_{j}x^{2}_{1})[D^*]$ is connected, we therefore have $wx^{1}_{j} \in E(G)$. Since $G$ is connected, the vertex $x^{2}_{1}$ is adjacent to a vertex, say $v$, that belongs to $A$. Since the vertex $x^{2}_{1}$ is adjacent neither $w$ nor $y$, the vertices $v$, $w$ and $y$ are distinct.

By Claim~\ref{lem p3} with $R = \{v,w\}$, there exists a path $P_{R,y}$ from $c_{0}$ to $y$ containing every vertex in $A_{0} \setminus \{v,w\}$. The path $P_{R,y}$ can be extended to a hamiltonian path of $G$ by proceeding along the edge $yx^{1}_{j-1}$ from $y$ to $x^{1}_{j-1}$, and then following the $(x^{1}_{j-1},x^{1}_{j})$-path, say $P^*$, on $C$ that does not contain the edge $x^{1}_{j-1}x^{1}_{j}$ (and contains all vertices of $P^1$) from $x^{1}_{j-1}$ to $x^{1}_{j}$, and then following the path $x^{1}_{j} w v x^2_1$ from $x^{1}_{j}$ to $x^2_1$; that is, the path
\[
c_0 P_{R,y} y, x^{1}_{j-1} P^* x^1_j, w v x^{2}_{1}
\]
is a hamiltonian path in $G$.  This completes the proof of Claim~\ref{lem z201}.3.~\smallqed

\medskip
By Claim~\ref{lem z201}.3, we may assume that $|V(P^{2})| \ge 2$, for otherwise $G$ is traceable and the desired result holds. Recall that $D_{1,2} \cap \barA = \{y\}$ and $D_{1,2} \cap \barA = \{x^{1}_{1}\}$. Further, the vertex $x^{1}_{1}$ is adjacent to neither $x^{2}_{1}$ nor $x^{2}_{n_2}$. In particular, $x^{1}_{1}x^{2}_{n_{2}} \notin E(G)$. In order for the set $D_{1,2}$ to dominate the vertex $x^{2}_{n_2}$, we note that $yx^{2}_{n_{2}} \in E(G)$.

We now consider the graph $G + x^{1}_{1}x^{2}_{n_{2}}$. For notational simplicity, let $D^* = D_{x^{1}_{1}x^{2}_{n_{2}}}$. By Claim~\ref{lem p2}, we have $|D^* \cap \barA| = |D^* \cap A| = 1$. Further, $|D^* \cap \{x^{1}_{1}, x^{2}_{n_{2}}\}| = 1$. Let $D^* \cap A = \{u\}$. By Lemma~\ref{lem 1}(c), the vertex $u$ is adjacent to  exactly one of $x^{1}_{1}$ and $x^{2}_{n_{2}}$. Therefore since $y$ is adjacent to both $x^{1}_{1}$ and $x^{2}_{n_{2}}$, we note that $u \ne y$.

Suppose firstly that $x^{1}_{1} \in D^*$. In this case, $u$ is adjacent to $x^{1}_{1}$ but not to $x^{2}_{n_{2}}$. Since $x^{1}_{1}x^{2}_{1} \notin E(G)$, in order for the set $D^*$ to dominate the vertex $x^{2}_{1}$, we note that $ux^{2}_{1} \in E(G)$. By Claim~\ref{lem p3} with $R = \{u\}$, there exists a path $P_{R,y}$ from $c_{0}$ to $y$ containing every vertex in $A_{0} \setminus \{u\}$. The path $P_{R,y}$ can be extended to a hamiltonian path of $G$ by proceeding along the edge $yx^{2}_{n_2}$ from $y$ to $x^{2}_{n_2}$, following the path $P^2$ in reverse direction from $x^{2}_{n_2}$ to $x^{2}_{1}$, proceeding along the path $x^{2}_{1}ux^1_1$ from $x^{2}_{1}$ to $x^1_1$, and then following the path $P^1$ from $x^1_1$ to $x^1_{n_1}$; that is, the path
\[
c_0 P_{R,y} y, x^{2}_{n_2}P^2x^{2}_{1}, u, x^{1}_{1} P^1 x^{1}_{n_1}
\]
is a hamiltonian path in $G$. Suppose next that $x^{2}_{n_{2}} \in D^*$. In this case, $u$ is adjacent to $x^{2}_{n_{2}}$ but not to $x^{1}_{1}$. Since $x^{1}_{n_1}x^{2}_{n_{2}} \notin E(G)$, in order for the set $D^*$ to dominate the vertex $x^{1}_{n_1}$, we note that $ux^{1}_{n_1} \in E(G)$. By Claim~\ref{lem p3} with $R = \{u\}$, there exists a path $P_{R,y}$ from $c_{0}$ to $y$ containing every vertex in $A_{0} \setminus \{u\}$. The path $P_{R,y}$ can be extended to a hamiltonian path of $G$ by proceeding along the edge $yx^{1}_{1}$ from $y$ to $x^{1}_{1}$, following the path $P^1$ from $x^1_1$ to $x^1_{n_1}$, proceeding along the path $x^1_{n_1} u x^2_{n_2}$ from $x^1_{n_1}$ to $x^1_{n_2}$, and following the path $P^2$ in reverse direction from $x^{2}_{n_2}$ to $x^{2}_{1}$; that is, the path
\[
c_0 P_{R,y} y, x^{1}_{1}P^1x^{1}_{n_1}, u, x^{2}_{n_2} P^2 x^{2}_{1}
\]
is a hamiltonian path in $G$. This completes the proof of Claim~\ref{lem z201}.~\smallqed

\medskip
By Claim~\ref{lem z201}, we may assume that $x^{2}_{1} \in D_{1,2}$, for otherwise $G$ is traceable and the desired result follows. Thus in this case, $D_{1,2} \cap (A \cup \barA) = \{x^{2}_{1},y\}$. Since $(G + x^{1}_{1}x^{2}_{1})[D_{1,2}]$ is a connected graph, $x^{2}_{1}y \in E(G)$. Moreover by Lemma~\ref{lem 1}(c), $yx^{1}_{1} \notin E(G)$. Since $x^{1}_{n_1}x^{2}_{1} \notin E(G)$, in order for the set $D_{1,2}$ to dominate the vertex $x^{1}_{n_1}$, we note that $yx^{1}_{n_1} \in E(G)$.

\begin{subclaim}
\label{lem z202}
If $|V(P^{2})| = 1$, then $G$ is traceable.
\end{subclaim}
\proof Suppose that $|V(P^{2})| = 1$. Thus, $P^{2}$ consists of the single vertex $x^{2}_{1}$. We show firstly that $x^{1}_{1}x^{1}_{n_{1}} \notin E(G)$. Suppose, to the contrary, that $x^{1}_{1}x^{1}_{n_{1}} \in E(G)$. We show that $x^{2}_{1}$ is not adjacent to any vertex of $P^{1}$. By our earlier observations, the vertex $x^{2}_{1}$ is adjacent to neither $x^1_1$ nor $x^1_{n_1}$. Suppose that $x^{2}_{1}$ is adjacent to $x^1_j$ for some $j$ where $1 < j < n_1$. In this case, $|V(P^{1})| \ge 3$, and so $G[\, \barA \,]$ has a cycle $C \colon P^{1} + x^{1}_{1}x^{1}_{n_{1}}$ as a subgraph. The $(x^{1}_{j+1},x^{1}_{j})$-path on $C$ that does not contain the edge $x^{1}_{j}x^{1}_{j+1}$ can be extended to a longer path in $G[\, \barA \,]$ by adding to it the vertex $x^{2}_{1}$ and the edge $x^{1}_{j}x^{2}_{1}$, contradicting the maximality of the path $P^1$. Hence, the vertex $x^{2}_{1}$ is adjacent to no vertex of $P^{1}$. Thus since $\{x^{2}_{1},y\}$ dominates all vertices of $P^1$ different from $x^1_1$, this implies that $y \succ P^{1} - x^{1}_{1}$. Therefore, $\{y, x^{1}_{n_{1}}\} \succ_{c} H$, and so $\gamma_c(H) \le 2$, a contradiction. Hence, $x^{1}_{1}x^{1}_{n_{1}} \notin E(G)$.

We now consider the graph $G + x^{1}_{n_{1}}x^{2}_{1}$. For notational simplicity, let $D^* = D_{x^{1}_{n_{1}}x^{2}_{1}}$. By Claim~\ref{lem p2}, we have $|D^* \cap \barA| = |D^* \cap A| = 1$. Further, $|D^* \cap \{x^{1}_{n_{1}}, x^{2}_{1}\}| = 1$. Let $D^* \cap A = \{u\}$. By Lemma~\ref{lem 1}(c), the vertex $u$ is adjacent to  exactly one of $x^{1}_{n_1}$ and $x^{2}_{1}$. Therefore since $y$ is adjacent to both $x^{1}_{n_1}$ and $x^{2}_{1}$, we note that $u \ne y$. If $x^{2}_{1} \in D^*$, then since $x^{2}_{1}x^{1}_{1} \notin E(G)$, we have $ux^{1}_{1} \in E(G)$. If $x^{1}_{n_{1}} \in D^*$, then since $x^{1}_{1}x^{1}_{n_{1}} \notin E(G)$, we must have $ux^{1}_{1} \in E(G)$. In both cases, $ux^{1}_{1} \in E(G)$.

By Claim~\ref{lem p3} with $R = \{y\}$, there exists a path $P_{R,u}$ from $c_{0}$ to $u$ containing every vertex in $A_{0} \setminus \{y\}$. The path $P_{R,u}$ can be extended to a hamiltonian path of $G$ by proceeding along the edge $ux^{1}_{1}$ from $u$ to $x^{1}_{1}$, following the path $P^1$ from $x^1_1$ to $x^1_{n_1}$, and then proceeding along the path $x^1_{n_1} y x^2_{1}$ from $x^1_{n_1}$ to $x^2_{1}$; that is, the path
\[
c_0 P_{R,u} u, x^{1}_{1}P^1x^{1}_{n_1}, y, x^{2}_{1}
\]
is a hamiltonian path in $G$. This completes the proof of Claim~\ref{lem z202}.~\smallqed

\medskip
By Claim~\ref{lem z202}, we may assume that $|V(P^{2})| \ge 2$, for otherwise $G$ is traceable and the desired result holds. We now consider the graph $G + x^{1}_{n_{1}}x^{2}_{1}$. For notational simplicity, let $D^* = D_{x^{1}_{n_{1}}x^{2}_{1}}$. By Claim~\ref{lem p2}, we have $|D^* \cap \barA| = |D^* \cap A| = 1$. Further, $|D^* \cap \{x^{1}_{n_{1}}, x^{2}_{1}\}| = 1$. Let $D^* \cap A = \{u\}$. By Lemma~\ref{lem 1}(c), the vertex $u$ is adjacent to  exactly one of $x^{1}_{n_1}$ and $x^{2}_{1}$. Therefore since $y$ is adjacent to both $x^{1}_{n_1}$ and $x^{2}_{1}$, we note that $u \ne y$.

Suppose firstly that $x^{2}_{1} \in D^*$. In this case, $u$ is adjacent to $x^{2}_{1}$ but not to  $x^{1}_{n_1}$. Since $x^{1}_{1}x^{2}_{1} \notin E(G)$, in order for the set $D^*$ to dominate the vertex $x^{1}_{1}$, we note that $ux^{1}_{1} \in E(G)$. By Claim~\ref{lem p3} with $R = \{u\}$, there exists a path $P_{R,y}$ from $c_{0}$ to $y$ containing every vertex in $A_{0} \setminus \{u\}$. The path $P_{R,y}$ can be extended to a hamiltonian path of $G$ by proceeding along the edge $yx^{1}_{n_1}$ from $y$ to $x^{1}_{n_1}$, following the path $P^1$ in the reverse direction from $x^1_{n_1}$ to $x^1_1$, proceeding along the path $x^1_{1} u x^2_{1}$ from $x^1_{1}$ to $x^2_{1}$, and following the path $P^2$ from $x^{2}_{1}$ to $x^{2}_{n_2}$; that is, the path
\[
c_0 P_{R,y} y, x^{1}_{n_1}P^1x^{1}_{1}, u, x^{2}_{1} P^2 x^{2}_{n_2}
\]
is a hamiltonian path in $G$. Suppose next that $x^{1}_{n_1} \in D^*$. In this case, $u$ is adjacent to $x^{1}_{n_1}$ but not to $x^{2}_{1}$. Since $x^{1}_{n_1}x^{1}_{n_2} \notin E(G)$ and $x^{2}_{1} \ne x^{2}_{n_{2}}$, in order for the set $D^*$ to dominate the vertex $x^{2}_{n_2}$, we note that $ux^{2}_{n_2} \in E(G)$.

By Claim~\ref{lem p3} with $R = \{y\}$, there exists a path $P_{R,u}$ from $c_{0}$ to $u$ containing every vertex in $A_{0} \setminus \{y\}$. The path $P_{R,u}$ can be extended to a hamiltonian path of $G$ by proceeding along the edge $ux^{2}_{n_2}$ from $u$ to $x^{2}_{n_2}$, following the path $P^2$ in the reverse direction from $x^2_{n_2}$ to $x^2_1$, proceeding along the path $x^2_{1} y x^1_{n_1}$ from $x^2_{1}$ to $x^1_{n_1}$, and following the path $P^1$ in the reverse direction from $x^{1}_{n_1}$ to $x^{1}_{1}$; that is, the path
\[
c_0 P_{R,u} u, x^{2}_{n_2}P^2x^{2}_{1}, y, x^{1}_{n_1} P^1 x^{1}_{1}
\]
is a hamiltonian path in $G$. This completes the proof of Claim~\ref{cl 2}.~\smallqed

\medskip
By Claims~\ref{cl 1} and~\ref{cl 2}, we may assume that $z \ge 3$, for otherwise $G$ is traceable and the desired result holds. The following claim uses similar ideas to those presented in~\cite{KaAn10}. However for completeness, we provide a proof of this claim.

\begin{claim}
\label{lem z30}
If $I$ is an independent set of $\barA$ where $|I| = t \ge 3$, then all the vertices of $I$ can be ordered $u_{1}, u_{2}, \ldots, u_{t}$ in such a way that there exist $t - 1$ different vertices $v_{1}, v_{2}, \ldots, v_{t - 1}$ of $A$ satisfying
$\{u_{i}, v_{i}\} \succ_{c} H - u_{i + 1}$
for all $i \in [t-1]$.
\end{claim}
\proof We will construct a tournament $T$ (a digraph which any two vertices are joined by an arc) with vertex set $V(T) = I$ and where the arcs of $T$ are defined as follow. For every two distinct vertices $u$ and $v$ in $I$, we choose a fixed $\gamma_{c}$-set, say $D_{uv}$, of $G + uv$. By Claim~\ref{lem p2}, $|D_{uv} \cap (A \cup \barA)| = 2$, implying that $|D_{uv} \cap \{u,v\}| = 1$ and $|D_{uv} \cap A| = 1$. Let $D_{uv} \cap A = \{x\}$. If $u \in D_{uv}$, then since $A$ is a complete subgraph, it follows that $\{u, x\} \succ_{c} H - v$. In this case, we orient the arc from $u$ to $v$. If $v \in D_{uv}$, then $\{v, x\} \succ_{c} H - u$, and we orient the arc from $v$ to $u$. We do this for every two distinct vertices $u$ and $v$ in $I$. This defines the arcs of the resulting tournament $T$. Since every tournament has a directed hamiltonian path, we let $u_{1} u_{2} \ldots u_{t}$ be a directed hamiltonian path in $T$. This implies that there exists a vertex $v_{i} \in A$ such that $\{u_{i}, v_{i}\} \succ_{c} H - u_{i + 1}$ for every $i \in [t-1]$. Since $I = \{u_{1},u_2,\ldots,u_{t}\}$ is an independent set, it follows that the vertex $v_i$ is adjacent to every vertex in $I$ except for the vertex $u_{i+1}$ for all $i \in [t-1]$. This implies that the vertices $v_1,v_2, \ldots, v_t$ are all distinct.~\smallqed

\medskip
We now return to the proof of Theorem~\ref{thm traceability}. Since $\{x^{1}_{1}, x^{2}_{1}, \ldots, x^{z}_{1}\}$ is an independent set of size~$z \ge 3$, by Claim~\ref{lem z30} there exists an ordering $u_{1}, u_{2}, \ldots , u_{z}$ of the vertices of $\{x^{1}_{1}, x^{2}_{1}, \ldots, x^{z}_{1}\}$ such that there exist vertices $v_{1}, v_{2}, \ldots, v_{z - 1}$ of $A$ satisfying
\[
\{u_{i}, v_{i}\} \succ_{c} H - u_{i + 1}
\]
for all $i \in [z-1]$. Let $R = \{v_{1}, v_{2}, \ldots, v_{z - 1}\}$. For notational convenience, if $u_{j} = x^{i}_{1}$ for some $i, j \in [z]$, then we relabel the path $P^{i}$ as the path $T^{j}$. Further we let $u'_{j} = x^{i}_{n_{i}}$. We note that $v_{i} \succ \{u_{1}, u_{2}, \ldots, u_{z}\} \setminus \{u_{i+1}\}$ for all $i \in [z-1]$. Further, we note that the collection of paths $P^1, P^2, \ldots, P^z$ is therefore precisely the collection of paths $T^1, T^2, \ldots, T^z$. If $|V(T^{i})| = 1$ for all $i \in [z]$, then $\{v_{1}, v_{2}\} \succ_{c} H$, and so $\gamma_c(H) \le 2$, a contradiction. Hence, $|V(T^{i})| > 1$ for some $i \in [z]$. For $i \in [z]$, we let
\[
\ell = \max \{i \, \colon |V(T^{i})| > 1 \}.
\]
Therefore if $\ell < z$, then $|V(T^{j})| = 1$ for all $j$ where $\ell < j \le z$, implying that $u_{j} = u'_{j}$ for such values of $j$. We remark that it is possible that $|V(T^{j})| = 1$ for some $j < \ell$. The structure of $G[ \, \barA \, ]$ is now illustrated by Figure~\ref{f:fig7}.

\setlength{\unitlength}{0.8cm}
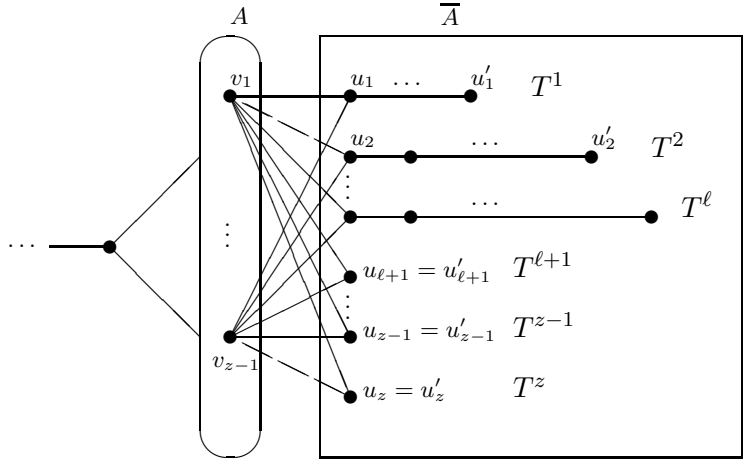
\begin{figure}[htb]
\begin{center}
\begin{picture}(13, 7)
\put(1.5, 2.5){\circle*{0.2}}
\put(1.5, 2.5){\line(1, 1){1.5}}
\put(1.5, 2.5){\line(1, -1){1.5}}
\put(3.5, 2.5){\oval(1, 7)}
\put(5, -1){\framebox(7,7){}}
\put(1.5, 2.5){\line(-1, 0){1}}
\put(3.5, 5){\circle*{0.2}}
\put(3.5, 1){\circle*{0.2}}
\put(3.4, 2.5){\footnotesize$\vdots$}
\put(3.5, 5.2){\footnotesize$v_{1}$}
\put(3.2, 0.5){\footnotesize$v_{z - 1}$}
\put(3.5, 5){\line(1, 0){2}}
\multiput(3.7, 4.9)(0.6, -0.3){3}{\line(2,-1){0.5}}
\put(3.5, 5){\line(1, -2){2}}
\put(3.5, 5){\line(1, -1){2}}
\put(3.5, 5){\line(2, -3){2}}
\put(3.5, 5){\line(2, -5){2}}
\put(3.5, 1){\line(1, 0){2}}
\multiput(3.7, 0.9)(0.6, -0.3){3}{\line(2,-1){0.5}}
\put(3.5, 1){\line(1, 1){2}}
\put(3.5, 1){\line(2, 1){2}}
\put(3.5, 1){\line(2, 3){2}}
\put(3.5, 1){\line(1, 2){2}}
\put(5.5, 5){\circle*{0.2}}
\put(5.5, 4){\circle*{0.2}}
\put(5.5, 3){\circle*{0.2}}
\put(6.5, 4){\circle*{0.2}}
\put(6.5, 3){\circle*{0.2}}
\put(5.5, 2){\circle*{0.2}}
\put(5.5, 1){\circle*{0.2}}
\put(5.5, 0){\circle*{0.2}}
\put(6.2, 5.2){\footnotesize$\ldots$}
\put(7.5, 4.2){\footnotesize$\ldots$}
\put(7.5, 3.2){\footnotesize$\ldots$}
\put(7.5, 5){\circle*{0.2}}
\put(9.5, 4){\circle*{0.2}}
\put(10.5, 3){\circle*{0.2}}
\put(5.5, 5){\line(1, 0){2}}
\put(5.5, 4){\line(1, 0){4}}
\put(5.5, 3){\line(1, 0){5}}
\put(5.5, 5.2){\footnotesize$u_{1}$}
\put(7.5, 5.2){\footnotesize$u'_{1}$}
\put(5.5, 4.2){\footnotesize$u_{2}$}
\put(9.5, 4.2){\footnotesize$u'_{2}$}
\put(5.4, 3.3){\footnotesize$\vdots$}
\put(5.4, 1.3){\footnotesize$\vdots$}
\put(-0.2, 2.5){\footnotesize$\ldots$}
\put(3.5, 6.2){\footnotesize$A$}
\put(7, 6.2){\footnotesize$\barA $}
\put(8.5, 5){$T^{1}$}
\put(10.5, 4){$T^{2}$}
\put(11, 3){$T^{\ell}$}
\put(8.2, 2){$T^{\ell + 1}$}
\put(5.7, 2){\footnotesize$u_{\ell + 1} = u'_{\ell + 1}$}
\put(8.2, 1){$T^{z - 1}$}
\put(5.7, 1){\footnotesize$u_{z - 1} = u'_{z - 1}$}
\put(8.2, 0){$T^{z}$}
\put(5.7, 0){\footnotesize$u_{z} = u'_{z}$}
\end{picture}
\vskip 0.75 cm
\caption{\label{f:fig7}
The structure of $G[ \, \barA \, ]$ after rearranging the paths}
\end{center}
\end{figure}

\begin{claim}
\label{cl l > 1}
If $\ell > 1$, then $G$ is traceable.
\end{claim}
\proof Suppose that $\ell > 1$. Therefore, $u'_{\ell} \ne u'_{1}$ and $u'_{\ell} \ne u_{\ell}$. We now consider the graph $G + u_{1}u'_{\ell}$. For notational simplicity, let $D_{1,\ell} = D_{u_{1}u'_{\ell}}$. By Claim~\ref{lem p2}, we have $|D_{1,\ell} \cap \barA| = |D_{1,\ell} \cap A| = 1$. Further, $|D_{1,\ell} \cap \{u_{1}, u'_{\ell}| = 1$. Let $D_{1,\ell} \cap A = \{w\}$. By Lemma~\ref{lem 1}(c), the vertex $w$ is adjacent to exactly one of $u_{1}$ and $u'_{\ell}$. By Claim~\ref{lem z30}, $\{u_{i}, v_{i}\} \succ_{c} H - u_{i + 1}$ for all $i \in [z-1]$. Therefore since $u_{i}u'_{i - 1} \notin E(G)$, this implies that $v_{i}u'_{i - 1} \in E(G)$ for all $i \in [z] \setminus \{1\}$.

Suppose firstly that $u_1 \in D_{1,\ell}$. Thus, $w$ is adjacent to $u_{1}$ but not to $u'_{\ell}$. If $w = v_{j}$ for some $j \in [z - 1]$, then $\{u_{1}, w\}$ does not dominate the vertex $u_{j + 1}$, contradicting the fact that $D_{1,\ell}$ is a (connected) dominating set of $G + u_1u'_{\ell}$. Hence, $w \notin R$. Since $u_{1}u_{z} \notin E(G)$, we note that $wu_{z} \in E(G)$. Further since $u_{1}u'_{z} \notin E(G)$ and $\{u_{1}, v_{1}\} \succ_{c} H - u_{2}$ and $z \ge 3$, it follows that $v_{1}u'_{z} \in E(G)$. (Possibly, $u_{z} = u'_{z}$.) As observed earlier, $v_{i}u'_{i - 1} \in E(G)$ for all $i \in [z] \setminus \{1\}$. By Claim~\ref{lem p3}, there exists a path $P_{R,w}$ from $c_{0}$ to $w$ containing every vertex in $A_{0} \setminus R$. The path $P_{R,w}$ can be extended to a hamiltonian path of $G$ by proceeding along the edge $wu_z$ from $w$ to $u_z$, following the path $T^z$ from $u_z$ to $u_z'$, proceeding along the path $u_z' v_1 u_1$ from $u_z'$ to $u_1$, following the path $T^1$ from $u_1$ to $u_1'$, proceeding along the path $u_1' v_2 u_2$ from $u_1'$ to $u_2$, following the path $T^2$ from $u_2$ to $u_2'$, proceeding along the path $u_2' v_3 u_3$ from $u_3$ to $u_3'$, and, continuing in this way, finally proceeding along the path $u_{z-2}' v_{z-1} u_{z-1}$ from $u_{z-2}'$ to $u_{z-1}$, and then following the path $T^{z-1}$ from $u_{z-1}$ to $u_{z-1}'$; that is, the path
\[
c_0 P_{R,w} u_z, u_z T^z u_z', v_1, u_1T^1u_1', v_2, u_2T^2u_2', \cdots, u_{z-2}T^{z-2}u_{z-2}', v_{z-1}, u_{z-1}T^{z-1}u_{z-1}'
\]
is a hamiltonian path in $G$, as desired. Suppose next that $u'_{\ell} \in D_{1,\ell}$. Thus, $w$ is adjacent to $u'_{\ell}$ but not to $u_{1}$. If $\ell \le z-2$, then $u'_{\ell}$ is adjacent to neither $u'_{z - 1}$ nor $u'_{z}$, implying that $w$ is adjacent to both $u'_{z - 1}$ and $u'_{z}$. If $\ell = z-1$, then $w$ is adjacent to $u'_{z-1}$. Further since $u'_{z-1}$ is not adjacent to $u'_{z}$, the vertex $w$ is adjacent to $u'_{z}$ in this case. If $\ell = z$, then $w$ is adjacent to $u'_{z}$. Further since $z \ge 3$, the vertex $u'_{z-1} \ne u_1$. Thus since the vertex $u'_{z}$ is not adjacent to $u'_{z-1}$, the vertex $w$ is therefore adjacent to $u'_{z-1}$ in this case. Thus in all cases, we note that the vertex $w$ is adjacent to both $u'_{z - 1}$ and $u'_{z}$.

Let $R_1 = (R \setminus \{v_1\}) \cup \{w\}$. By Claim~\ref{lem p3}, there exists a path $P_{R_1}$ from $c_{0}$ to $v_1$ containing every vertex in $A_{0} \setminus R_1$. The path $P_{R_1,v_1}$ can be extended to a hamiltonian path of $G$ by proceeding along the edge $v_1u_1$ from $v_1$ to $u_1$, following the path $T^1$ from $u_1$ to $u_1'$, proceeding along the path $u_1' v_2 u_2$ from $u_1'$ to $u_2$, following the path $T^2$ from $u_2$ to $u_2'$, proceeding along the path $u_2' v_3 u_3$ from $u_3$ to $u_3'$, and, continuing in this way, finally proceeding along the path $u_{z-2}' v_{z-1} u_{z-1}$ from $u_{z-2}'$ to $u_{z-1}$, and then following the path $T^{z-1}$ from $u_{z-1}$ to $u_{z-1}'$, proceeding along the path $u_{z-1}' w u_{z}'$ from $u_{z-1}'$ to $u_{z}'$, following the path $T^z$ in the reverse direction from $u_z'$ to $u_z$; that is, the path
\[
c_0 P_{R_1,v_1} v_1, u_1 T^1 u_1', v_2, u_2T^2u_2', \cdots, u_{z-1}T^{z-1}u_{z-1}', w, u_{z}'T^{z}u_{z}
\]
is a hamiltonian path in $G$, as desired.~\smallqed

\medskip
By Claim~\ref{cl l > 1}, we may assume that $\ell = 1$, for otherwise $G$ is traceable and the desired result follows. Thus, $P^{1} = T^{1}$ and $|V(P^{1})| = n_1 \ge 2$, and so $u_{1} \ne u'_{1}$. Moreover, $u_{i} = u'_{i}$ for all $i \in [z] \setminus \{1\}$.

Suppose firstly that $u_{1}u'_{1} \in E(G)$, and so  $G[\, \barA \,]$ contains a cycle $C \colon T^{1} + u_{1}u'_{1}$ as a subgraph. If there exist integers $j$ and $r$ where $j \in [n_1 - 1] \setminus \{1\}$ and $r \in [z] \setminus \{1\}$ such that $u_{r}x^{1}_{j} \in E(G)$, then the $(x^{1}_{j+1},x^{1}_{j})$-path on $C$ that does not contain the edge $x^{1}_{j}x^{1}_{j+1}$ can be extended to a longer path in $G[\, \barA \,]$ by adding to it the vertex $u_r$ and the edge $u_{r}x^{1}_{j}$, contradicting the maximality of the path $P^1$. Hence, no vertex of $T^{1}$ is adjacent to any vertex from the set $\{u_{2}, u_{3}, \ldots, u_{z}\}$, implying that the vertex $u_i$ is an isolated vertex in $G[\, \barA \,]$ for all $i \in [z] \setminus \{1\}$. Since $\{u_{2}, v_{2}\} \succ_{c} H - u_{3}$ and $u_2$ is isolated in $G[\, \barA \,]$, this implies that $v_{2} \succ H - u_{3}$. Therefore, $\{v_{1}, v_{2}\} \succ_{c} H$, and so $\gamma_c(H) \le 2$, a contradiction. Hence, $u_{1}u'_{1} \notin E(G)$.

We now consider the graph $G + u_{1}u'_{1}$. For notational simplicity, let $D_{1,1} = D_{u_{1}u'_{1}}$. By Claim~\ref{lem p2}, we have $|D_{1,1} \cap \barA| = |D_{1,1} \cap A| = 1$. Further, $|D_{1,1} \cap \{u_{1}, u'_{1}| = 1$. Let $D_{1,1} \cap A = \{w\}$. By Lemma~\ref{lem 1}(c), the vertex $w$ is adjacent to exactly one of $u_{1}$ and $u'_{1}$. If $w = v_{i}$ for some $i \in [z - 1]$, then $D_{1,1}$ does not dominate $u_{i + 1}$, a contradiction. Hence, $w \in A \setminus R$. Since neither $u_1$ nor $u_1'$ is adjacent to $u_z$ or $u_{z - 1}$, the vertex $w$ is necessarily adjacent to both  $u_z$ and $u_{z - 1}$. Recall that $v_{i}u_{i-1} \in E(G)$ for all $i \in [z - 1] \setminus \{1\}$, and recall that $v_{2}u'_{1} \in E(G)$. Let $R_1 = (R \setminus \{v_1\}) \cup \{w\}$. By Claim~\ref{lem p3}, there exists a path $P_{R_w}$ from $c_{0}$ to $v_1$ containing every vertex in $A_{0} \setminus R_1$. The path $P_{R_1,v_1}$ can be extended to a hamiltonian path of $G$ by proceeding along the edge $v_1u_1$ from $v_1$ to $u_1$, following the path $T^1$ from $u_1$ to $u_1'$, and then proceeding along the path $u_1' v_2 u_2 v_3 \ldots v_{z-1} u_{z-1} w u_z$; that is, the path
\[
c_0 P_{R_1,v_1} v_1, u_1 T^1 u_1' v_2 u_2 v_3 \cdots v_{z-1} u_{z-1} w  u_{z}
\]
is a hamiltonian path in $G$. This completes the proof Theorem~\ref{thm traceability}.~\qed

\section{$k$-$\gamma_{c}$-Critical Graphs which are Non-Traceable}

In this section, we establish the realizability result that for $k \ge 4$, there exist $k$-$\gamma_{c}$-critical graphs which is non-traceable containing $\zeta$ vertices for all $0 \le \zeta \le k - 4$. For this purpose, for $k \ge 3$ we introduce a class $\cP(k)$ of $k$-$\gamma_{c}$-critical graphs such that, for every graph $G \in \cP(k)$ and every integer $\ell \ge 1$, there exists a $(k + \ell)$-$\gamma_{c}$-critical graph that contains $G$ as an induced subgraph. Further, we construct a class $\cN(s)$ of graphs for all $s \ge 6$.

\noindent
\textbf{The class $\cP(k)$ for $k \ge 3$.} A $k$-$\gamma_{c}$-critical graph $G$ is in the class $\cP(k)$ if there exists a maximal complete subgraph $H$ of $G$ of order at least~$2$ satisfies the following two properties.
\\ [-26pt]
\begin{enumerate}
\item Every vertex of $G$ belongs to some $\gamma_{c}$-set of $G$ that contains a vertex of $H$.
\item For every pair of non-adjacent vertices $x$ and $y$ in $G$, there exists a CD-set $D'_{xy}$ of $G + xy$ such that $D_{xy}' \cap V(H) \ne \emptyset$ and $|D'_{xy}| < k$ (we remark that $D'_{xy}$ need not necessarily be a $\gamma_{c}$-set of $G + xy$).
\end{enumerate}


\noindent \textbf{The class $\cN(s)$ for $s \ge 6$}. For a set $S = [s]$ where $s \ge 6$, we let $B_{1} = \{a_{i} \colon i \in [s]\}$ and $B_{2} = \{b_{i} \colon i \in [s]\}$ be two disjoint sets of vertices, and let
\[
B_{3} = \left\{z_{i,j} \colon \{i, j\} \in \dbinom{S}{2}\right\}
\]
where $\dbinom{S}{2}$ is a set of all pairs (regardless of order) of the members in $S$, and so $|B_{3}| = \binom{s}{2}$. A graph $G$ in the class $\cN(s)$ can be constructed from the disjoint sets $B_{1}$, $B_{2}$ and $B_{3}$ by adding a new vertex $x$ and adding edges as follows: \\ [-22pt]
\begin{enumerate}
\item[$\bullet$] Add edges so that $B_{1}$ and $B_{2}$ form two complete subgraphs.
\item[$\bullet$] Add all edges between $B_{1}$ and $B_{2}$ except for the edges $a_ib_i$ for $i \in [s]$.
\item[$\bullet$] Join $x$ to every vertex of $B_{3}$.
\item[$\bullet$] Join $b_{i}$ to $z_{j, \ell}$ for $1 \le i \ne j \ne \ell \le s$.
\item[$\bullet$] Join $a_{i}$ to $z_{i,j}$ for $1 \le i \ne j \le s$.
\end{enumerate}
We note that for $i \in [s]$, $N_{B_{3}}(a_{i}) = \{z_{i, 1}, z_{i, 2}, \ldots, z_{i, i - 1}, z_{i, i + 1}, z_{i, i + 2}, \ldots, z_{i, s}\}$ and
\[
N_{B_{3}}(b_{i}) = \left\{z_{j, \ell} : \{j, \ell\} \in \dbinom{S \setminus \{i\}}{2}\right\}.
\]
A graph in the class $\cN(s)$ is illustrated by Figure~\ref{f:fig8}.

\setlength{\unitlength}{0.8cm}
\begin{figure}[htb]
\begin{center}
\begin{picture}(9, 7)
\put(4, 6){\circle*{0.2}}
\put(0, 4){\circle*{0.2}}
\put(2, 4){\circle*{0.2}}
\put(4, 4){\circle*{0.2}}
\put(8, 4){\circle*{0.2}}
\put(2, 2){\circle*{0.2}}
\put(2, 0){\circle*{0.2}}

\put(4, 6){\line(-2, -1){4}}
\put(4, 6){\line(-1, -1){2}}
\put(4, 6){\line(0, -1){2}}
\put(4, 6){\line(2, -1){4}}
\put(2, 0){\line(1, 2){2}}
\put(2, 0){\line(-1, 2){2}}
\put(2, 2){\line(3, 1){6}}

\multiput(2, 0)(0.6, 0.4){10}{\line(3, 2){0.5}}

\qbezier(2, 0)(0.8, 2)(2, 4)

\put(4, 2){\oval(5, 1)}
\put(4, 0){\oval(5, 1)}

\put(6, 4){\footnotesize$\ldots$}
\put(2.8, 4){\footnotesize$\ldots$}
\put(3.9, 2){\footnotesize$\ldots$}
\put(3.9, 0){\footnotesize$\ldots$}
\put(2.2, 1.8){\footnotesize$b_{1}$}
\put(2.2, -0.2){\footnotesize$a_{1}$}

\put(-0.8, 4.2){\footnotesize$z_{1, 2}$}
\put(1.2, 4.2){\footnotesize$z_{1, 3}$}
\put(4.2, 4.2){\footnotesize$z_{1, s}$}
\put(8, 4.2){\footnotesize$z_{s - 1, s}$}

\put(4, 6.2){\footnotesize$x$}

\put(7, 0){\footnotesize$B_{1}$}
\put(7, 2){\footnotesize$B_{2}$}
\put(9.5, 4){\footnotesize$B_{3}$}

\multiput(2, 0)(0, 0.6){7}{\line(0,1){0.5}}
\multiput(2.1, 2.1)(0.6, 0.6){3}{\line(1,1){0.5}}
\multiput(1.9, 2.1)(-0.6, 0.6){3}{\line(-1,1){0.5}}
\end{picture}
\vskip 0.5 cm
\caption{\label{f:fig8}
A graph $G$ in the class $\cN(s)$}
\end{center}
\end{figure}

Let $\cF(k, \zeta)$ be the class of $k$-$\gamma_{c}$-critical graphs with $\zeta$ cut-vertices which are non-traceable. In view of Theorem \ref{thm hamiltonian}, $\cF(k, \zeta) = \emptyset$ for all $k \in [3]$. We show next that $\cN(s) \subseteq \cP(4)$ and $\cN(s) \subseteq \cF(4, 0)$. In particular, this implies that the class $\cP(4)$ is not empty when $k = 4$.

\begin{lem}\label{lem NT1}
For all $s \ge 6$, $\cN(s) \subseteq \cF(4, 0)$. Moreover, $\cN(s) \subseteq \cP(4)$ where in the construction of $\cP(4)$ here we take $H$ as the maximal complete subgraph $G[B_{2}]$.
\end{lem}
\proof Let $G \in \cN(s)$. We show that $G$ is a $4$-$\gamma_{c}$-critical non-traceable graph. Let $H$ be the maximal complete subgraph $G[B_{2}]$ of $G$. We show firstly that $\gamma_c(G) \ge 4$. Suppose, to the contrary, that there exists a CD-set $D$ of $G$ of size~$3$. Suppose that $x \in D$. If $D = \{x, z_{i, j}, z_{i', j'}\}$, then $D$ does not dominate $a_{\ell}$ where $\ell \in S \setminus \{i, j, i', j'\}$. If $D = \{x, z_{i, j}, a_{i}\}$, then $D$ does not dominate $b_{i}$. If $D = \{x, z_{i, j}, b_{\ell}\}$, then $D$ does not dominate $a_{\ell}$. In all three cases we produce a contradiction. Hence, $x \notin D$. In order to dominate the vertex~$x$, we have $z_{i, j} \in D$ for some $i$ and $j$ where $1 \le i \ne j \le s$. If $D = \{z_{i, j}, a_{i}, a_{j}\}$, then $D$ does not dominate $z_{i', j'}$ where $\{i', j'\} \cap \{i, j\} = \emptyset$. If $D = \{z_{i, j}, a_{i}, a_{\ell}\}$ where $\ell \notin \{i,j\}$, then $D$ does not dominate $z_{i', j'}$ where $\{i', j'\} \cap \{i, j, \ell\} = \emptyset$. If $D = \{z_{i, j}, a_{i}, b_{\ell}\}$ where $\ell \notin \{i,j\}$, then $D$ does not dominate $z_{j, \ell}$. If $D = \{z_{i, j}, a_{i}, z_{i, j'}\}$ or $D = \{z_{i, j}, b_{\ell}, z_{i', j'}\}$, then $D$ does not dominate $z_{j, \ell}$ where $\ell \notin \{i,i',j,j'\}$. If $D = \{z_{i, j}, b_{\ell}, b_{\ell'}\}$ where $\ell, \ell' \notin \{i,j\}$, then $D$ does not dominate $z_{\ell,\ell'}$. In all cases, we have a contradiction. We deduce, therefore, that not such CD-set $D$ of size~$3$ exists. Hence, $\gamma_c(G) \ge 4$.

We show next that property~(a) holds in the construction of $\cP(4)$, and, simultaneously, we show that $\gamma_{c}(G) = 4$. For all $i$ and $j$ where $1 \le i \ne j \le s$, we note that $x z_{i, j} a_{i} b_{j}$ is an induced path in $G$ and the set $D_{i,j} = \{x, z_{i, j}, a_{i}, b_{j}\}$ is a CD-set of $G$. This implies that $\gamma_c(G) \le 4$ and every vertex of $G$ belongs to some CD-set of $G$ of size~$4$ that contains a vertex of $H$, where recall that $V(H) = B_2$. As observed earlier, $\gamma_c(G) \ge 4$. Consequently, $\gamma_c(G) = 4$ and every vertex of $G$ belongs to some $\gamma_{c}$-set of $G$ that contains a vertex of $H$. This establishes property~(a) in the construction of $\cP(4)$.

We show next that property~(b) in the construction of $\cP(4)$ holds. Let $u$ and $v$ be an arbitrary pair of non-adjacent vertices of $G$. We show that there exists a CD-set $D_{uv}$ of $G + uv$ such that $|D_{uv}| = 3$ and the set $D_{uv}$ contains at least one vertex in $B_{2}$.

Suppose that $x \in \{u, v\}$. Renaming vertices if necessary, we may assume that $x = u$. Thus, $v \in B_{1} \cup B_{2}$. If $v = a_{i}$, then let $D_{uv} = \{x, a_{i}, b_{j}\}$ where $i \ne j$. If $v = b_{i}$, then let $D_{uv} = \{x, a_{i}, b_{j}\}$ where $i \ne j$. In both cases, $|D_{uv}| = 3$, the set $D_{uv}$ contains a vertex of $B_2$ and $D_{uv} \succ_{c} G + uv$, as desired. Hence, we may assume that $x \notin \{u, v\}$, for otherwise the desired result holds.

Suppose next that $z_{i, j} \in \{u, v\}$ for some $i$ and $j$ where $1 \le i \ne j \le s$. Renaming vertices if necessary, we may assume that $z_{i, j} = u$. If $v = a_{\ell}$, then necessarily $\ell \ne i$ and we let $D_{uv} = \{a_{\ell}, b_{\ell}, z_{i, j}\}$. If $v = b_{i}$ or $v = b_j$, then we let $D_{uv} =  \{z_{i, j}, b_{i}, b_{j}\}$. If $v = z_{i', j'}$ where $\{i, j\} \ne \{i', j'\}$, then $z_{i, j}$ is adjacent to at least one of the vertices $b_{i'}$ or $b_{j'}$ and we let $D_{uv} =  \{z_{i, j}, b_{i'}, b_{j'}\}$. In all cases, $|D_{uv}| = 3$, the set $D_{uv}$ contains a vertex of $B_2$ and $D_{uv} \succ_{c} G + uv$, as desired. Hence, we may assume that $\{u, v\} \subseteq B_{1} \cup B_{2}$, for otherwise the desired result holds. Thus, $\{u, v\} = \{a_{i}, b_{i}\}$ for some $i \in [s]$. We now let $D_{uv} = \{a_{i}, b_i, z_{i, j}\}$ where $i \ne j$. Once again in this case, $|D_{uv}| = 3$, the set $D_{uv}$ contains a vertex of $B_2$ and $D_{uv} \succ_{c} G + uv$, as desired.

Thus, for an arbitrary pair $x$ and $y$ of non-adjacent vertices of $G$, there exists a CD-set $D_{uv}$ of $G + uv$ such that $|D_{uv}| = 3$ and the set $D_{uv}$ contains at least one vertex of $H$, where $V(H) = B_2$. As observed earlier, $\gamma_c(G) = 4$. Therefore, property~(b) in the construction of $\cP(4)$ holds. In particular, we note $G$ is a $4$-$\gamma_{c}$-critical graph. Since $G$ is an arbitrary graph in $\cN(s)$, we have $\cN(s) \subseteq \cP(4)$.

By construction, we note that $G$ has no cut-vertex. Finally, to show that $\cN(s) \subseteq \cF(4, 0)$, it remains to show that $G$ is non-traceable. Let $S = B_{1} \cup B_{2} \cup \{x\}$ and consider the graph $G - S$. We note that $|S| = 2s+1$ and $G - S$ consists of $\binom{s}{2}$ isolated vertices, namely the vertices $z_{i, j}$ where $1 \le i \ne j \le s$. Since $s \ge 6$, we note that $2s + 2 < \binom{s}{2} = \frac{s(s - 1)}{2}$. Hence,
\[
|S| = 2s + 1 < 2s + 2 < \frac{s(s + 1)}{2} = |B_{3}| = \omega(G - S).
\]
Therefore, by Observation~\ref{o:toughoftrace} the graph $G$ is non-traceable. Thus, $\cN(s) \subseteq \cF(4, 0)$. This completes the proof of Lemma~\ref{lem NT1}.~\qed

\medskip
For $k \ge 4$ and for $\ell \ge 1$, we next give a construction of a $(k + \ell)$-$\gamma_{c}$-critical graph that contains a graph $G$ in the class $\cP(k)$ as an induced subgraph. Let $G$ be a graph in the class $\cP(k)$, and let $H$ be a maximal complete subgraph of $G$ having properties (a) and (b) in the construction of the class $\cP(k)$. For $\ell \ge 1$, let $G_1, \ldots, G_{\ell}$ be $\ell$ vertex disjoint complete graphs where $G_i = K_{n_{i}}$ and $n_{i} \ge 1$ for $i \in [\ell]$. Let $G(n_{1}, n_{2}, \ldots, n_{\ell})$ be the graph constructed from an isolated vertex $x_{0}$, vertex disjoint copies of the complete graphs $G_1, \ldots, G_{\ell}$ and $G \in \cP(k)$ by adding edges according to the join operations
\[
x_{0} \vee G_1 \vee G_2 \vee \cdots \vee G_\ell \vee \leftidx{_H}{G}.
\]
Let $\cP(k,\ell)$ be the class of all such graphs $G(n_{1}, n_{2}, \ldots, n_{\ell})$. A graph in the class $\cP(k,\ell)$ is illustrated in Figure~\ref{f:fig9}.

\setlength{\unitlength}{0.8cm}
\begin{figure}[htb]
\begin{center}
\begin{picture}(17.25, 5)
\put(3.8, 3){\circle*{0.2}}
\put(3.8, 3){\line(1, 1){0.7}}
\put(3.8, 3){\line(1, -1){0.7}}
\put(9, 3){\oval(1, 2)}
\put(9.6, 2.5){\line(1, 0){0.7}}
\put(9.6, 3.5){\line(1, 0){0.7}}
\put(9.6, 3){\line(1, 0){0.7}}
\put(5, 3){\oval(1, 2)}
\put(5.6, 2.5){\line(1, 0){0.7}}
\put(5.6, 3.5){\line(1, 0){0.7}}
\put(5.6, 3){\line(1, 0){0.7}}
\put(7, 3){\oval(1, 2)}
\put(11, 3){\oval(1, 2)}
\put(12, 3){\oval(3, 3)}
\put(3.7, 2.2){\footnotesize$x_{0}$}
\put(4.7, 2.9){\footnotesize$K_{n_{1}}$}
\put(6.7, 2.9){\footnotesize$K_{n_{2}}$}
\put(8.7, 2.9){\footnotesize$K_{n_{\ell}}$}
\put(4.8, 1.5){\footnotesize$G_1$}
\put(6.8, 1.5){\footnotesize$G_2$}
\put(8.8, 1.5){\footnotesize$G_\ell$}
\put(7.8, 2.95){\footnotesize$\ldots$}
\put(11.3, 1){\footnotesize$G \in \cP(k)$}
\put(10.8, 2.9){\footnotesize$H$}
\end{picture}
\vskip -0.85 cm
\caption{\label{f:fig9}
The graph in the class $\cP(k,\ell)$}
\end{center}
\end{figure}
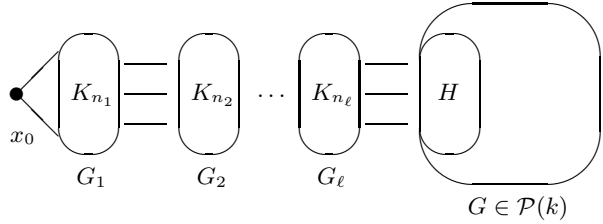


\begin{thm}\label{thm gl}
For $k \ge 4$ and $\ell \ge 1$, every graph in the class $\cP(k,\ell)$ is a $(k + \ell)$-$\gamma_{c}$-critical graph.
\end{thm}
\proof For $k \ge 4$ and $\ell \ge 1$, let $G(n_{1}, n_{2}, \ldots, n_{\ell})$ be a graph in the class $\cP(k,\ell)$ that is constructed from a graph $G \in \cP(k)$ with $H$ as the maximal complete subgraph of $G$ having properties (a) and (b) in the construction of $G$. For notation convenience, we write the graph $G(n_{1}, n_{2}, \ldots, n_{\ell})$ simply as $G_{k,\ell}$. Let $G_1, \ldots, G_{\ell}$ be the $\ell$ vertex disjoint complete graphs used to construct $G_{k,\ell}$, where $G_i = K_{n_{i}}$ for $i \in [\ell]$. Let $x_i$ be an arbitrary vertex of $G_i$ for $i \in [\ell]$, and let $X = \{x_1,\ldots,x_\ell\}$. Let $D$ be a $\gamma_{c}$-set of $G$ that contains a vertex of $H$. Thus, $|D| = \gamma_c(G) = k$ and $D \cap V(H) \ne \emptyset$.

We show firstly that $\gamma_{c}(G_{k,\ell}) = k + \ell$. Since $D \cup X \succ_{c} G$, we note that $\gamma_{c}(G_{k,\ell}) \le |D| + |X| = k + \ell$. To show that $\gamma_{c}(G_{k,\ell}) \ge k + \ell$, let $D'$ be an arbitrary $\gamma_{c}$-set of $G_{k,\ell}$. If $D'$ contains the vertex $x_0$, then since $D'$ is a CD-set of $G_{k,\ell}$ it also contains a vertex of $G_1$. If $D'$ does not contains the vertex $x_0$, then in order to dominate the vertex $x_{0}$, we note that $D'$ contains a vertex of $G_1$. Hence in both cases, $D' \cap V(G_1) \ne \emptyset$. Since $G - H$ is not the empty graph, the set $D'$ contain at least one vertex of $G$. By the connectedness of $G_{k,\ell}[D']$, the set $D'$ therefore contains at least one vertex from each of the sets $G_i$ for $i \in [\ell]$ and $D' \cap V(H) \ne \emptyset$. Therefore, $|D' \cap (\cup^{\ell}_{i = 1}V(G_i))| \ge \ell$, that is, $|D' \cap (V(G_{k,\ell}) \setminus V(G))| \ge \ell$. Since $H$ is a complete subgraph of $G$ and $D' \cap V(H) \ne \emptyset$, we note that the set $D' \cap V(G)$ is a CD-set of $G$, implying that $|D' \cap V(G)| \ge \gamma_c(G) = k$. Hence, $\gamma_c(G_{k,\ell}) = |D'| = |D' \cap V(G)| + |D' \cap (V(G_{k,\ell}) \setminus V(G))| \ge k + \ell$. Consequently, $\gamma_{c}(G_{k,\ell}) = k + \ell$.

We establish next the criticality of $G_{k,\ell}$. Let $u$ and $v$ be an arbitrary pair of non-adjacent vertices of $G_{k,\ell}$. We show that there exists a CD-set $D_{uv}$ of $G_{k,\ell} + uv$ such that $|D_{uv}| < k + \ell$.
We first consider the case when $\{u, v\} \subseteq V(G_{k,\ell}) \setminus V(G)$. Renaming vertices if necessary, we may assume that $u = x_{j}$ and $v = x_{j'}$ where $0 \le j < j' \le \ell$. Since $u$ and $v$ are non-adjacent vertices of $G_{k,\ell}$, we note that $j + 2 \le j'$. If $j > 0$, then let $D_{uv} = D \cup (X \setminus \{x_{j+1}\})$. If $j = 0$, then let $D_{uv} = D \cup (X \setminus \{x_{1}\})$. In both cases, $|D_{uv}| = k + \ell - 1$ and $D_{uv} \succ_{c} G_{k,\ell} + uv$. Hence, we may assume that at least one of $u$ and $v$ belongs to $G$, for otherwise the desired result follows. Renaming vertices if necessary, we may assume that $u \in V(G)$. By property~(a) in the construction of the graph $G \in \cP(4)$, there exists a $\gamma_{c}$-set $D_u$ of $G$ that contains the vertex~$u$ and contains a vertex of $H$.

Suppose next that $v \notin V(G)$, and so $v = x_0$ or $v \in V(G_i)$ for some $i \in [\ell]$. Renaming vertices if necessary, we may assume that $v = x_i$ for some $i \in [\ell] \cup \{0\}$.
Suppose that $u \in V(H)$, implying that $v \ne x_{\ell}$.  If $v = x_{0}$, then let $D_{uv} = D_u \cup (X \setminus \{x_{1}\})$. If $v \ne x_{0}$, then let $D_{uv} = D_u \cup (X \setminus \{x_{\ell}\})$. In both cases, $|D_{uv}| = k + \ell - 1$ and $D_{uv} \succ_{c} G_{k,\ell} + uv$. Hence, we may assume that $u \in V(G) \setminus V(H)$. Once again if $v \ne x_{0}$, then let $D_{uv} = D_u \cup (X \setminus \{x_{\ell}\})$, and if $v = x_{i}$ for some $i \in [\ell]$, then let $D_{uv} = D_u \cup (X \setminus \{x_{\ell}\})$. In both cases, $|D_{uv}| = k + \ell - 1$ and $D_{uv} \succ_{c} G_{k,\ell} + uv$. Hence, we may assume that $v \in V(G)$. By property~(b) in the construction of the graph $G \in \cP(4)$, there exists a $\gamma_{c}$-set $D_{uv}'$ of $G$ that contains a vertex of $H$ and such that $|D_{uv}'| \leq k - 1$. In this case, we let $D_{uv} = D_{uv}' \cup X$ and note that $|D_{uv}| \leq |D_{uv}'| + |X| \leq k + \ell - 1$ and $D_{uv} \succ_{c} G_{k,\ell} + uv$. These observations imply that $G_{k,\ell}$ is a $(k + \ell)$-$\gamma_{c}$-critical graph. This completes the proof of Theorem~\ref{thm gl}.~\qed

\medskip
We are now ready to establish the realisability of $k$-$\gamma_{c}$-critical non-traceable graphs containing $\zeta$ cut-vertices for all $k \ge 4$ and $0 \le \zeta \le k - 4$. Recall that $\cF(k, \zeta)$ is the class of $k$-$\gamma_{c}$-critical graphs with $\zeta$ cut-vertices which are non-traceable.

\begin{thm}\label{thm real}
For integers $k \ge 1$ and $\zeta \ge 0$, $\cF(k, \zeta) \ne \emptyset$ if and only if $k \ge 4$ and $0 \le \zeta \le k - 4$.
\end{thm}
\proof We first show that if $k \ge 4$ and $0 \le \zeta \le k - 4$, then $\cF(k, \zeta) \ne \emptyset$. In view of Lemma~\ref{lem NT1}, $\cF(4, 0) \ne \emptyset$. Thus if $k = 4$ and $\zeta = k - 4 = 0$, then $\cF(k, \zeta) \ne \emptyset$. Hence, we may assume that $k \ge 5$, for otherwise the desired result follows. Let $G \in \cN(s)$ for some integer $s \ge 6$. Adopting our earlier notation, Lemma~\ref{lem NT1} yields also that $G \in \cP(4)$ where here we take $H$ as the maximal complete subgraph $G[B_{2}]$ having properties (a) and (b) in the construction of~$G$.
For a given $\zeta \in [k - 4] \cup \{0\}$, let $G^* = G(n_{1}, n_{2}, \ldots, n_{k-4})$ be a graph in the class $\cP(4,k-4)$ that is constructed from the graph $G$ by taking $n_{i} \ge 2$ for $i \in [k-4]$ in the case when $\zeta = 0$ and taking $n_{i} = 1$ for $i \in [\zeta]$ and $n_{i} \ge 2$ for $i \in [k-4] \setminus [\zeta]$ in the case when $\zeta \in [k - 4]$ for the $k - 4$ complete graphs $K_{n_{1}}, K_{n_{2}}, \ldots, K_{n_{k - 4}}$ used in the construction of $G^*$. We note that if $\zeta = 0$, then $G^*$ has no cut-vertex, while if $\zeta \in [k - 4]$, then $G^*$ has exactly $\zeta$ cut-vertices, namely the singleton vertices of the $\zeta$ complete graphs $K_{n_{1}}, \ldots, K_{n_{\zeta}}$. In both cases, $G^*$ has exactly $\zeta$ cut-vertices. Moreover, by Theorem~\ref{thm gl} the graph $G^*$ is a $k$-$\gamma_{c}$-critical.

We show next that $G^*$ is non-traceable. Adopting our earlier notation used in the construction of the graph $G \in \cP(4)$, let $S = B_{1} \cup B_{2} \cup \{x\}$ and consider the graph $G^* - S$. We note that $|S| = 2s+1$ and $G^* - S$ consists of $\binom{s}{2}$ isolated vertices, namely the vertices $z_{i, j}$ where $1 \le i \ne j \le s$, together with an additional component containing the vertex $x_0$ and the $k - 4$ complete graphs $K_{n_{1}}, K_{n_{2}}, \ldots, K_{n_{k - 4}}$. Since $s \ge 6$, we note that $2s + 2 < \binom{s}{2} = \frac{s(s - 1)}{2}$. Hence,
\[
|S| = 2s + 1 < 2s + 2 < \frac{s(s + 1)}{2} + 1 = |B_{3}| + 1  = \omega(G^* - S).
\]

Therefore, by Observation~\ref{o:toughoftrace} the graph $G^*$ is non-traceable, implying that $G^* \in \cF(k, \zeta)$. Hence, if $k \ge 4$ and $0 \le \zeta \le k - 4$, then $\cF(k, \zeta) \ne \emptyset$.

Conversely, suppose that $\cF(k, \zeta) \ne \emptyset$. Since every hamiltonian graph is traceable, by Theorem~\ref{thm hamiltonian}, we must have that $k \ge 4$. Theorem~\ref{thm mpm} implies that $\cF(k, \zeta) = \emptyset$ when $\zeta \ge k - 1$. Thus, $\zeta \le k - 2$. By Theorem~\ref{thm traceability} also implies that $\cF(k, \zeta) = \emptyset$ when $\zeta \in \{k-3,k-2\}$. These results imply that $0 \le \zeta \le k - 4$ and $k \ge 4$. This completes the proof of Theorem~\ref{thm real}.~\qed


\medskip
\begin{thebibliography}{99}
\label{bib}\small

\bibitem{An07} N. Ananchuen, On domination critical graphs with cut vertices having connected domination number~$3$. \textit{Int. Math. Forum} \textbf{2} (2007), 3041--3052.

\bibitem{AAP} N. Ananchuen, W. Ananchuen, and M. D. Plummer, Matching properties in connected domination critical graphs. \textit{Discrete Math.} \textbf{308} (2008), 1260--1267.

\bibitem{BiFe14} D. Binkele-Raible and H. Fernau, A parameterized measure-and-conquer analysis for finding a k-leaf spanning tree in an undirected graph. \textit{Discrete Math. Theor. Comput. Sci.} \textbf{16}(1) (2014), 179--200.

\bibitem{CaWeYu00} Y. Caro, D. B. West, and R. Yuster, Connected domination and spanning trees with many leaves. \textit{SIAM J. Discrete Math.} \textbf{13}(2) (2000), 202--211.

\bibitem{ChSuMa04} X. G. Chen, L. Sun, and D. X. Ma, Connected domination critical graphs. \textit{Applied Math. Letters} \textbf{17} (2004), 503--507.


\bibitem{Ch73} V. Chv\'{a}tal, Tough graphs and hamiltonian circuits. \textit{Discrete Math.} \textbf{5} (1973), 215--228.

\bibitem{CoSt90} C. J. Colbourn and L. K. Stewart, Permutation graphs: connected domination and Steiner trees. \textit{Discrete Math.} \textbf{86} (1990), 179--189.

\bibitem{DaMo88} A. D'Atri and M. Moscarini, Distance-hereditary graphs, Steiner trees, and connected domination. \textit{SIAM J. Comput.} \textbf{17}(3) (1988), 521--538.

\bibitem{DeHaHe13} W. J. Desormeaux, T. W. Haynes, and M. A. Henning, Bounds on the connected domination number of a graph. \textit{Discrete Appl. Math.} \textbf{161}(18) (2013), 2925--2931.

\bibitem{DuMa09} W. Duckworth and B. Mans, Connected domination of regular graphs. \textit{Discrete Math.} \textbf{309}(8) (2009), 2305--2322.



\bibitem{FaTiZh97} O. Favaron, F. Tian and Lei. Zhang, Independence and hamiltonicity in 3-domination-critical graphs. \textit{J. Graph Theory} \textbf{25} (1997), 173--184.

\bibitem{Fe09} M. Fellows, D. Lokshtanov, N. Misra, M. Mnich, F. Rosamond, and S. Saurabh, The complexity ecology of parameters: an illustration using bounded max leaf number. \textit{Theory Comput. Systems} \textbf{45}(4) (2009), 822--848.

\bibitem{FlTiWeZh99} E. Flandrin, F. Tian, B. Wei and L. Zhang, Some properties of 3-domination-critical graphs. \textit{Discrete Math.} \textbf{205} (1999), 65--76.




\bibitem{Ga94} G. Galbiati, F. Maffioli, and A. Morzenti, A short note on the approximability of the maximum leaves spanning tree problem. \textit{Inf. Process. Lett.}   \textbf{52}(1) (1994), 45--49.

\bibitem{HHS1}  T. W.\ Haynes, S. T. Hedetniemi and P. J. Slater, \emph{Fundamentals of Domination in Graphs}, Marcel Dekker, Inc., New York, 1998.


\bibitem{HHS2} T. W. Haynes, S. T. Hedetniemi and P. J. Slater (eds), \emph{Domination  in Graphs: Advanced Topics}, Marcel Dekker, Inc. New York, 1998.



\bibitem{HeYe_book} M. A. Henning and A. Yeo, \emph{Total Domination in Graphs (Springer Monographs in Mathematics)} 2013. ISBN: 978-1-4614-6524-9 (Print) 978-1-4614-6525-6 (Online).

\bibitem{thesis} P. Kaemawichanurat, Connected domination critical graphs. Ph.D. Thesis. PhD supervisor: Louis Caccetta. Curtin University (2016).



\bibitem{Ka19} P. Kaemawichanurat, Connected domination critical graphs with $k - 3$ cut vertices, submitted.

\bibitem{KaAn10} P. Kaemawichanurat and N. Ananchuen, On $4$-$\gamma_{c}$-critical graphs with cut vertices. \textit{Utilitas Math.} \textbf{82} (2010), 253--268.


\bibitem{KaAn19} P. Kaemawichanurat and N. Ananchuen, Connected domination critical graphs with cut vertices. To appear in \textit{Discuss. Math. Graph Theory.}


\bibitem{KaCa19} P. Kaemawichanurat and L. Caccetta, Hamiltonicity of domination critical claw-free graphs. To appear in \textit{J. Combin. Comput. Combin. Math.}


\bibitem{KaCaAn18} P. Kaemawichanurat, L. Caccetta and W. Ananchuen, Hamiltonicity of connected domination critical graphs. \textit{Ars Combin.} \textbf{136} (2018), 127--151.

\bibitem{KJ}  P. Kaemawichanurat and T. Jiarasuksakun, Some results on the independence number of connected domination critical graphs. AKCE Int. J. Graphs Comb. 15 (2018), no. 2, 190--196.

\bibitem{Ka12} H. Karami, S. M. Sheikholeslami, A. Khodkar, and D. B. West, Connected domination number of a graph and its complement. \textit{Graphs Combin.} \textbf{28}(1) (2012), 123--131.

\bibitem{LiZh15} T. Liu, Z. Lu, and K. Xu, Tractable connected domination for restricted bipartite graphs. \textit{J. Comb. Optim.} \textbf{29}(1) (2015), 247--256.

\bibitem{SaWa79} E. Sampathkumar and H. B. Walikar, The connected domination number of a graph. \textit{J. Math. Phys. Sci.} \textbf{13}(6) (1979), 607--613.

\bibitem{Sa91} L. A. Sanchis, Maximum number of edges in connected graphs with a given domination number. \textit{Discrete Math.} \textbf{8}(1) (1991), 65--72.


\bibitem{SuBl83} D. P. Sumner and P. Blitch, Domination critical graphs. \textit{J. COmbin. Theory B} \textbf{34} (1983) 65--76.

\bibitem{TiWeZh99} F. Tian, B. Wei and L. Zhang, Hamiltonicity in $3$-domination critical graphs with $\alpha = \delta + 2$. \textit{Discrete Appl. Math.} \textbf{92} (1999), 57--70.


\bibitem{Wo90} E. Wojcicka, Hamiltonian properties of domination critical graphs. \textit{J. Graph Theory} \textbf{14} (1990), 205--215.

\bibitem{WuLi99} J. Wu and H. Li, On calculating connected dominating set for efficient routing in ad hoc wireless networks. \emph{Proceedings of the 3rd International Workshop on Discrete Algorithms and Methods for Mobile Computing and Communications} ACM (1999), 7--14.

\bibitem{YuChXiYoXi05} Y. Yuansheng, Z. Chengye, L. Xiaohui, J. Yongsong and H. Xin, Some $3$-connected $4$-edge critical non-Hamiltonian graphs\textit{J. Graph Theory} \textbf{50} (2005), 316--320.

\end{thebibliography}
\end{document}